\documentclass[12pt,reqno]{amsart}
\usepackage{amsfonts}
\usepackage{amssymb}
\usepackage{amsthm}
\usepackage{upref}
\usepackage{enumerate}

\usepackage{amscd}

\makeatletter
\def\LaTeX{\leavevmode L\raise.42ex
    \hbox{\kern-.3em\size{\sf@size}{0pt}\selectfont A}\kern-.15em\TeX}
 \makeatother
 
\DeclareMathOperator{\clos}{clos}

\numberwithin{equation}{section}

\newtheorem{lemma}{Lemma}[section]
\newtheorem{theorem}[lemma]{Theorem} 
\newtheorem{corollary}[lemma]{Corollary}
\newtheorem{proposition}[lemma]{Proposition}

\theoremstyle{definition}

\newtheorem{example}[lemma]{Example}

\newtheorem{remark}[lemma]{Remark}

\newcommand{\arccosh}{\operatorname{arccosh}}

  \newcommand{\e}{\eqref}

\newcommand{\q}{\quad}

\newcommand{\ov}{\overline}
\newcommand{\wt}{\widetilde}

\newcommand{\ti}{\tilde}

\renewcommand{\d}{\delta}

   \newcommand{\sgn}{\operatorname{sgn}}

\renewcommand\Im{\operatorname{Im}}
\renewcommand\Re{\operatorname{Re}}

\newenvironment{pf}{\begin{proof}}{\end{proof}}

\def\qqq{\mathrel{\subset\mkern-15mu\lower.38ex\hbox{${\scriptscriptstyle\rightarrow}$}}}

\let\cal\mathcal

\let\Bbb\mathbb
 
\begin{document}
\title 
[Asymptotic  behavior of orthogonal polynomials]
{Asymptotic  behavior of orthogonal polynomials.   Singular
critical case} 
\author{ D. R. Yafaev  }
\address{   Univ  Rennes, CNRS, IRMAR-UMR 6625, F-35000
    Rennes, France and SPGU, Univ. Nab. 7/9, Saint Petersburg, 199034 Russia}
\email{yafaev@univ-rennes1.fr}
 

\thanks {Supported by  project   Russian Science Foundation   17-11-01126}

\begin{abstract}
Our goal  is to find  an asymptotic behavior as $n\to\infty$ of the orthogonal polynomials $P_{n}(z)$  defined by  Jacobi recurrence coefficients $a_{n}$ (off-diagonal terms) and $ b_{n}$ (diagonal terms).
We consider the case  $a_{n}\to\infty$, $b_{n}\to\infty$ in such a way that $\sum  a_{n}^{-1}<\infty$  
$($that is, the Carleman condition is violated$)$ and
$\gamma_{n}:=2^{-1}b_{n} (a_{n}a_{n-1})^{-1/2} \to \gamma $ as $n\to\infty$. In the case $|\gamma  | \neq 1$
asymptotic formulas for
$P_{n}(z)$ are known; they depend crucially on the sign of $| \gamma  |-1$. We study  the critical case $| \gamma  |=1$.
The formulas obtained are qualitatively different in the cases $|\gamma_{n}|
 \to  1-0$  and $|\gamma_{n}|
 \to  1+0$.  Another goal of the paper is to advocate an approach to a study of asymptotic behavior of  $P_{n}(z)$  based on a close analogy of the Jacobi difference equations  and differential equations of Schr\"odinger type.
      \end{abstract}



 \maketitle

\section{Introduction}


\subsection{Orthogonal polynomials and Jacobi operators}

 Orthogonal polynomials  $P_{n } (z)$ can be defined by  a recurrence relation
   \begin{equation}
 a_{n-1} P_{n-1} (z) +b_{n} P_{n } (z) + a_{n} P_{n+1} (z)= z P_n (z), \q n\in{\Bbb Z}_{+}, \q z\in{\Bbb C}, 
\label{eq:RR}\end{equation}
with the boundary conditions $P_{-1 } (z) =0$, $P_0 (z) =1$. We always suppose that $a_n >0$, $b_n=\bar{b}_n $. Determining $P_{n } (z)$, $n=1,2, \ldots$, successively from \e{eq:RR}, we see that $P_{n } (z)$ is a polynomial of degree $n$: $P_{n } (z)= p_{n}z^n+\cdots$ where $p_{n}= (a_{0}a_{1}\cdots a_{n-1})^{-1}$. 

 For all $z$ with $\Im z \neq 0$, the equation \e{eq:RR} either has exactly one (up to a multiplicative constant) solution in $  \ell^2 ({\Bbb Z}_{+})$ or all its solutions are in $  \ell^2 ({\Bbb Z}_{+})$. The first instance is known as the limit point case   and the second one  -- as  the  the limit circle case.

It is natural (see, e.g., the book \cite{AKH}) to associate with the coefficients $a_{n} , b_{n} $  
 a  three-diagonal matrix  
  \begin{equation}
{\cal J} = 
\begin{pmatrix}
 b_{0}&a_{0}& 0&0&0&\cdots \\
 a_{0}&b_{1}&a_{1}&0&0&\cdots \\
  0&a_{1}&b_{2}&a_{2}&0&\cdots \\
  0&0&a_{2}&b_{3}&a_{3}&\cdots \\
  \vdots&\vdots&\vdots&\ddots&\ddots&\ddots
\end{pmatrix} 
\label{eq:ZP+}\end{equation}
known as the Jacobi matrix.
Then equation \e{eq:RR} with the boundary condition $P_{-1 } (z) =0$ is
 equivalent to the equation ${\cal J} P( z )= z P(z)$ for the vector $P(z)=(P_{0}(z), P_1(z),\ldots)^{\top}$.  Thus $P(z)$ is an ``eigenvector" of the matrix $\cal J$ corresponding to an ``eigenvalue" $z$. 
  

Let us now consider Jacobi operators  defined   by matrix \e{eq:ZP+} in the canonical basis $e_{0}, e_{1}, \ldots$ of the space $\ell^2 ({\Bbb Z}_{+})$. 
The minimal Jacobi operator $J_{0}$ is  defined   by the formula $J_{0} u= {\cal J} u$ on a set $\cal D \subset \ell^2 ({\Bbb Z}_{+})$ of vectors $u$ with only a  finite number of non-zero components.  It is symmetric  in the space $\ell^2 ({\Bbb Z}_{+})$, and its adjoint operator $J^*_{0}$ is given by the same  formula $J_{0}^* u= {\cal J} u$ on all vectors $u \in \ell^2 ({\Bbb Z}_{+})$ such that ${\cal J}u \in \ell^2 ({\Bbb Z}_{+})$.  The operator $J_{0}$ is essentially self-adjoint in the limit point case, and  it has deficiency indices  $(1,1)$  in the limit circle case. In the limit point case, the closure $\clos J_{0}$ of $J_{0}$ and its adjoint operator are defined on the same domain: ${\cal D} (\clos J_{0})= {\cal D} (J_{0}^*)$.


   The spectra of all  self-adjoint extensions $J$ of the minimal operator $J_{0}$ are simple with $e_{0} = (1,0,0,\ldots)^{\top}$ being a generating vector. Therefore it is natural to define   the   spectral measure of $J$ by the relation $d\rho_{J}(\lambda)=d(E_{J}(\lambda)e_{0}, e_{0})$ where  $E_{J}(\lambda)$      is the spectral family of the operator $J$.
  For all   extensions $J$ of the operator $J_{0}$, the  polynomials  $P_{n}(\lambda)$ are orthogonal and normalized  in the spaces $L^2 ({\Bbb R};d\rho_{J})$: 
      \[
\int_{-\infty}^\infty P_{n}(\lambda) P_{m}(\lambda) d\rho_{J}(\lambda) =\d_{n,m};
\]
as usual, $\d_{n,n}=1$ and $\d_{n,m}=0$ for $n\neq m$.
 
 
 The comprehensive presentation of the results  described shortly above  can be found in  the books \cite{AKH, Chihara,   Schm} and the surveys \cite{Lub,Simon, Tot}.
     
   
      \subsection{Asymptotic behavior   of orthogonal polynomials}
   
We are interested in the asymptotic behavior of the polynomials  $P_{n } (z)$ as $n\to\infty$. It is of course to be expected that asymptotic formulas for $P_{n } (z)$ depend crucially on the behavior 
of recurrence coefficients $a_{n}$ and $b_{n}$  for large $n$.  A study of this problem was initiated by P.~Nevai in his book \cite{Nev}. He (see also the papers \cite{Mate} by A.~M\'at\'e,  P.~Nevai, and V.~Totik and  \cite{Va-As}  by W.~Van Assche and J.~S.~Geronimo) investigated the case where
$a_{n} \to a_{\infty}>0$, $b_{n} \to 0$   as $n\to\infty$.   

The case of the coefficients  $a_{n}\to \infty$     was later  studied in \cite{Jan-Nab} by J.~Janas and  S.~Naboko and in \cite{Apt} by A.~I.~Aptekarev and J.~S.~Geronimo. It was assumed in these papers  that the growth of  $a_{n}$ is not too rapid. More precisely, the so called Carleman condition 
    \begin{equation}
\sum_{n=0}^\infty a_{n}^{-1}=\infty
\label{eq:Carl}\end{equation}
(introduced by T.~Carleman in his book \cite{Carleman}) was required.
Under this assumption  for arbitrary  $b_{n}$,
the operators $J_{0}$ are essentially
 self-adjoint on $\cal D$. 
  With respect to the coefficients $b_{n}$, it was assumed in \cite{Jan-Nab, Apt}    that there exists a   limit
 \begin{equation}
   \frac{b_{n}}{2\sqrt{a_{n-1}a_{n}} }=: \gamma_{n}\to  \gamma ,  \q n\to\infty,
\label{eq:Gr}\end{equation}
 where $| \gamma | < 1$ so that $b_{n}$  are relatively small compared to $a_{n}$.  A typical result of these papers for a particular case $b_{n} =0$ is stated  below as Theorem~\ref{DD0}      (it was required in \cite{Apt} that $b_{n}\to\infty$, but this is probably an oversight).
 
 The case of rapidly increasing coefficients $a_{n}$ when the Carleman condition \e{eq:Carl} is violated,  so that
  \begin{equation}
  \sum_{n =0}^{\infty} a_{n}^{-1} <\infty ,
 \label{eq:nc}\end{equation} 
 was investigated in a recent paper \cite{nCarl}. 
 It was assumed in  \cite{nCarl} that $|\gamma | \neq 1$. The asymptotic formulas for $P_{n}  (z)$ turn out to be qualitatively different for $|\gamma | < 1$ and $|\gamma | > 1$.   Astonishingly, the asymptotics of the orthogonal polynomials in this a priori highly singular case is particularly simple and general.

  Let us briefly describe the results of \cite{nCarl}.
 In the case $|\gamma |<1$ the orthogonal polynomials are oscillating for large $n$: 
 \begin{equation}
P_{n}(  z )= a_{n} ^{-1/2}   \Big(\kappa_{+} (z)  e^{-i\phi_{n} }+ \kappa_{-}  (z) e^{ i\phi_{n}}+ o( 1)\Big), \q n\to\infty,
\label{eq:A2P}\end{equation} 
where
 \begin{equation}
 \phi_{n}= \sum_{m=0}^{n-1} \arccos \, (-\gamma_{m}), \q n\geq 1,
\label{eq:Grf}\end{equation}
and    $\kappa_{\pm} (z) \in{\Bbb C}$  are  some constants.  The same (only the constants $\kappa_{\pm} (z)$ change)  formula is true for all solutions of equation \e{eq:RR}. All of them
   belong to   $  \ell^2 ( {\Bbb Z}_{+})$  because the factor $\{a_{n} ^{-1/2}\}  \in  \ell^2 ( {\Bbb Z}_{+})$   due to    the condition \e{eq:nc}.  This implies that in the case $|\gamma |<1$, the Jacobi operators $J_{0}$ have deficiency indices $(1,1)$.
   
   In the case $|\gamma |>1$, the operators $J_{0}$ are essentially self-adjoint (provided the coefficients $a_{n}$ are polynomially bounded), and  the orthogonal polynomials are exponentially growing:
 \begin{equation}
    P_{n}(  z )= \kappa (z) a_{n} ^{-1/2} (- \sgn \gamma )^n e^{\varphi_{n}}  (1+ o(1)), \q n\to\infty, 
\label{eq:A2P1}\end{equation} 
(unless $z$ is an eigenvalue of the self-adjoint operator $J =\clos J_{0}$)  where  
 \begin{equation}
\varphi_{n}= \sum_{m=0}^{n-1} \arccosh |\gamma_{m}|, \q n\geq 1.
\label{eq:Grf1}\end{equation}
 We emphasize that a  finite number of terms in \e{eq:Grf} and \e{eq:Grf1} are arbitrary; only the values of $a_{n}, b_{n}$  for large $n$ are essential.

According to \e{eq:A2P} and \e{eq:A2P1} the asymptotic behavior of the polynomials $P_{n} (z)$
is the same for all $z\in {\Bbb C}$, both for real $z$ and for $z$ with $\Im z\neq 0$.   Only the coefficients $\kappa_{\pm}  (z) $ and $ \kappa  (z) $ depend on $z$.


 We also note   the paper \cite{Sw-Tr} where the Carleman and non-Carleman cases were treated   at an equal footing so that the difference in the corresponding asymptotic formulas for $P_{n } (z)$ was not quite visible.
 
 \section{Main results}
 
 The paper has two goals. The first one is to present a general approach to a study of asymptotic behavior of orthogonal polynomials based on an analogy with the theory of differential equations of Schr\"odinger type. The second goal is to apply this scheme to a study of the critical case where 
 $|\gamma | =1$ in \e{eq:Gr}. We concentrate here on rapidly increasing coefficients when
 condition \e{eq:nc} is satisfied.
 
 \subsection{Critical case}

  
  A classical example where the critical case $|\gamma | = 1$ occurs is given by  the Laguerre polynomials $L_{n}^{(p)}(z) $; those are  the orthogonal polynomials determined by the recurrence coefficients  
     \begin{equation}
    a_{n} = \sqrt{(n+1)(n+1+p)} \q\mbox{and}\q     b_{n} =  2n+p+1, \q p>-1. 
\label{eq:Lag}\end{equation}
The corresponding Jacobi operators $J=J^{(p)}$ have absolutely continuous  spectra
coinciding with $[0,\infty)$.

For sufficiently general coefficients $a_{n}, b_{n}$, the critical case was studied in the papers  \cite{Jan-Nab-Sh, Na-Si} (see also the references therein) where the Carleman condition  \e{eq:Carl} was required.
Our goal is to study the critical case for rapidly growing coefficients $a_{n}, b_{n}$ when the Carleman condition is not satisfied. The asymptotic formulas we obtain are quite different from those of the papers \cite{Jan-Nab-Sh, Na-Si}.

To handle the critical case, we require more specific  assumptions on the coefficients $a_{n}$ and $b_{n}$.  
To make our presentation as simple as possible, we assume that
  \begin{equation}
  a_{n}=  n^\sigma (1+  \alpha n^{-1}+ O (n^{-2}) )
\label{eq:ASa}\end{equation}
and
  \begin{equation}
  b_{n}=  2 \gamma n^\sigma  (1+  \beta n^{-1}+ O (n^{-2}) )
\label{eq:ASb}\end{equation}
for some   $\alpha,\beta, \gamma\in{\Bbb R}$. The critical case is distinguished  by the condition $|\gamma|=1$. We set
\[
\nu= -\sgn \gamma.
\]
As discussed in Sect.~1.2, in the non-critical case  $|\gamma|\neq 1$ asymptotic formulas are qualitatively different for  $\sigma<1$ and for $\sigma>1$. In the  critical case the borderline is
  $\sigma=3/2$. Here we study the singular situation $\sigma>3/2$.

 Note   that    if $P_{n}  (z)$ are the orthogonal polynomials corresponding to  coefficients $( a_{n}, b_{n})$, then according to equation \e{eq:RR} the polynomials $ (-1)^n P_{n}  (-z)$  correspond to the coefficients $(a_{n}, -b_{n})$. Therefore without loss of generality, we may suppose that $\gamma  = 1$ in \e{eq:Gr}  or \e{eq:ASb}.

The results below crucially depend on the value of
 \begin{equation}
\tau= 2 \beta -2\alpha + \sigma .
\label{eq:BX3}\end{equation}
As an example, note that for the Jacobi coefficients \e{eq:Lag}, we have $\sigma=1$, $\alpha=1+ p/2$, $\beta= (1+p)/2$, so that $\tau =0$.
In the main bulk of the paper we suppose that
$\tau\neq 0$.
   The case $\tau=0$ (doubly critical) is discussed in the final Sect.~6.

    Let us state our main results.  We first distinguish  solutions of the Jacobi equation
     \begin{equation}
 a_{n-1} f_{n-1} (z) +b_{n} f_{n } (z) + a_{n} f_{n+1} (z)= z f_n (z), \q n\in{\Bbb Z}_{+}, \q z\in{\Bbb C}, 
\label{eq:Jy}\end{equation}
by their behavior for $n\to\infty$.

\begin{theorem}\label{GSSx}
       Let the assumptions  \e{eq:ASa},  \e{eq:ASb} with $| \gamma |=1$  and $\sigma>3/2$  be satisfied.  Set $\varrho =\min\{\sigma-3/2,1/2\}$. 
     For all   $z\in  \Bbb C$,    the equation \e{eq:Jy} has a solution $\{f_{n}( z )\}$ with asymptotics
    \begin{equation}
f_{n}(  z )  =  \nu^n n^{-\sigma /2  + 1/4} e^{ - 2 i \sqrt{| \tau | n}} \big(1 + O( n^{-\varrho})\big) , \q n\to \infty,
\label{eq:A22G+}\end{equation} 
if $\tau<0$, and with asymptotics
    \begin{equation}
f_{n}(  z )  =  \nu^n n^{-\sigma /2  + 1/4} e^{  -2  \sqrt{\tau n}} \big(1 + O( n^{-\varrho})\big) , \q n\to \infty,
\label{eq:A22G-}\end{equation} 
if  $\tau>0$.
 Asymptotic relations  \e{eq:A22G+} and  \e{eq:A22G-}  are  uniform in $z$ from compact subsets    of the complex plane $  \Bbb C $.
For all $n\in {\Bbb Z}_{+}$, the functions $f_{n}( z )$ are entire functions of $z\in  \Bbb C$ of minimal exponential type.
 \end{theorem}
 
 By analogy with differential equations, it is natural to use the term ``Jost solutions" for solutions $\{f_{n} (z)\}$ constructed in Theorem~\ref{GSSx}.
 Note that according to  \e{eq:A22G+} and  \e{eq:A22G-} the leading   terms of their asymptotics  do not depend on $z$. This is a   unusual phenomenon  specific for rapidly growing coefficients $a_{n}$.

 For orthogonal polynomials, we have the following result.
 
  \begin{theorem}\label{GPS}
    Let the assumptions of Theorem~\ref{GSSx} be satisfied.    Then, for all $z\in{\Bbb C}$,  the sequence of the orthogonal polynomials $P_{n}( z )$  has  asymptotics
  \begin{equation}
P_{n}(  z )= \nu^n n^{-\sigma /2  + 1/4}   \Big(\kappa_{+} (z) e^{-2 i \sqrt{| \tau | n}  }+ \kappa_{-} (z) e^{ 2 i \sqrt{| \tau | n} }+ O( n^{-\varrho})\Big), \q n\to\infty,
\label{eq:A22P+}\end{equation} 
if $\tau<0$,  and
   \begin{equation}
  P_{n}(z)= \kappa (z)  \nu^n  n^{-\sigma /2  + 1/4} e^{2\sqrt{\tau  n}}  (1+ O( n^{-\varrho})), \q n\to\infty,
\label{eq:GEGE+}\end{equation}
if $\tau>0$.   Here $\kappa_{\pm}  (z)$ and $\kappa   (z)$ are some complex constants.  Asymptotic relations  \e{eq:A22P+} and \e{eq:GEGE+} are  uniform in $z$ from compact subsets    of the complex plane $  \Bbb C $.
 \end{theorem}


Note that the right-hand sides of  \e{eq:A22P+} and \e{eq:GEGE+} depend on $z$ only through   asymptotic constants $\kappa_{\pm} (z)$ and $\kappa (z)$. 
These constants      can be expressed via the Wronskians  of the polynomial  $\{P_{n}  (z)\}_{n=-1}^\infty$ and the Jost   $\{f_{n} (z)\}_{n=-1}^\infty$  solutions of the Jacobi equation  \e{eq:Jy} (see Sect.~5.1). 
 Obviously, formulas \e{eq:A22P+}  and \e{eq:GEGE+}  play the role of formulas \e{eq:A2P}  and \e{eq:A2P1}  for the non-critical case.
  
 Spectral results are stated in the following assertion.

  \begin{theorem}\label{S-Adj}
    Let the assumptions of Theorem~\ref{GSSx} be satisfied. 
   
$1^0$ If    $\tau<0$, then
the minimal Jacobi  operator $J_{0}$ has deficiency indices $(1,1)$ so that the spectra of all its self-adjoint extensions are discrete. 

  $2^0$  If  $\tau > 0$, then the operator  $J_{0}$ is essentially self-adjoint and the spectrum of  its closure is semi-bounded from below and discrete. 
 \end{theorem}
 
 We emphasize that although for all $\tau\neq 0$ the spectra of Jacobi operators are discrete, the reasons for this are different in the cases $\tau< 0$ and $\tau> 0$. In the first case, the discreteness of the spectra follows from general results of N.~Nevanlinna \cite{Nevan} on Jacobi operators that are not essentially self-adjoint. In the second case, diagonal elements $b_{n}$ dominate in some sense
 off-diagonal elements $a_{n}$, but diagonal operators always have   discrete spectra. 
 
  Note that the critical situation studied here is morally similar to a threshold behavior of orthogonal polynomials for the case $a_{n}\to a_{\infty}>0$ , $ b_{n}\to 0$ as $n\to\infty$. For such coefficients,  the role of \e{eq:Gr} is played (see \cite{Nev,Mate,JLR}) by the relation 
   \[
 \lim_{n\to\infty}  \frac{b_{n}-\lambda}{2 a_{n}} =-\frac{ \lambda}{2 a_\infty} .
\]
Since the essential spectrum of the operator $J$ is now $[-2 a_\infty,2 a_\infty]$, the values $\lambda=\pm 2 a_\infty$ are the threshold values of the spectral parameter $\lambda$.  The parameter $-\lambda/ (2 a_\infty)$ plays the role of $\gamma $ so that the cases $|\gamma |<1$ (resp., $|\gamma |>1$)  correspond to $\lambda$ lying inside the essential spectrum of   $J$  (resp., outside of it).
 
   \subsection{Classification}
 
 Generically,  leading terms of the Jost solutions asymptotics depend on the spectral parameter $z$. We call such situation ``generic" or ``regular". However, it might happen that the dependence on $z$ disappears in  asymptotic formulas for $f_{n}  (z)$. We call this situation ``exceptional" or ``singular".  It occurs if the coefficients $a_{n}\to\infty$ sufficiently rapidly. However the conditions on the growth of $a_{n}$ are different in non-critical and critical cases.
 
 Let us classify possible cases accepting assumptions \e{eq:ASa} and  \e{eq:ASb}.
 In  the non-critical case when   $|\gamma | \neq 1$    according to \e{eq:A2P} and \e{eq:A2P1},  the  situation  is singular if   $\sigma>1$. For the regular situation where $\sigma\leq 1$, we refer to the papers \cite{Jan-Nab, Apt}.

 In the critical case when $|\gamma | = 1$ and $\tau\neq 0$, we use Theorems~\ref{GSSx} or \ref{GPS}.
 Now  formulas  \e{eq:A22P+}  and \e{eq:GEGE+}  are true (and hence we are in the singular case) if   $\sigma>3/2$. The regular situation was studied  in  the papers  \cite{Jan-Nab-Sh, Na-Si}.
 
 In Sect.~6, we also consider the doubly critical case when $|\gamma | = 1$ and $\tau =0$.  Here the regular (singular) situation occurs if $\sigma\leq 2$ (resp., $\sigma > 2$).  The orthogonal polynomials with Jacobi coefficients  \e{eq:Lag}  fall into the regular case.
 

 

 
  
 The discussion above is  summarized in Figure~1. 
 
  
  \bigskip
  
 \begin{figure}  [hbt!]
 
 \centering
  
 \begin{tabular}{  l | l  | l  }

  \phantom & Regular & Singular
   \\
  \hline
  Non-critical:  $|\gamma | \neq 1$ & $\sigma \leq 1$   & $\sigma>1$ 
   \\
    \hline
Critical: $|\gamma | = 1$, $\tau \neq 0$    & $\sigma \leq 3/2$ & $\sigma>3/2$ 
   \\
    \hline
   Doubly critical: $|\gamma | = 1$, $\tau =0$ & $\sigma \leq 2$ & $\sigma>2$ 
   \\
    \hline

   \end{tabular} 

  \caption{Regular and singular  cases}
\end{figure}

 
 
  
  
   \subsection{Discrete versus continuous}
   
   
 Let us compare   difference  \e{eq:Jy} and  differential 
 \begin{equation}
 - (a(x) f ' (x, z) )'+ b(x) f  (x, z)= z f  (x, z), \q  x>0,  \q a(x) > 0,
\label{eq:Schr}\end{equation}
equations. To a large extent, $x$, $a(x)$ and $b(x)=\ov{b(x)}$ in \e{eq:Schr} play the roles of the parameters $n$, $a_n$ and $b_{n}$ in the Jacobi equation \e{eq:Jy}. The regular solution $\psi(x,z)$ of the differential equation \e{eq:Schr}  is distinguished by the conditions
 \[
\psi (0, z) =0,\q   \psi' (0, z) =1.
\]
It plays the role of the polynomial solution $P_n (z)$ of equation \e{eq:Jy} distinguished by the conditions
$P_{-1} (z)=0$, $P_0 (z)=1$.

A study of asymptotics of the regular solution $\psi(x,z)$  relies on a construction of   special solutions of the differential equation  \e{eq:Schr}   distinguished by their asymptotics as $x\to\infty$. For example, in the case $a(x)=1$, $b\in L^1 ({\Bbb R}_{+})$,  equation  \e{eq:Schr} has a solution $f  (x, z)$,  known as the Jost solution, behaving like $e^{i\sqrt{z}x}$, $  \Im \sqrt{z} \geq 0$, as $x\to\infty$.
Under fairly general assumptions  
equation \e{eq:Schr} has a solution $f(x,z)$ (we also  call it the Jost solution) whose asymptotics is given by  the classical Liouville-Green formula (see   Chapter~6 of the book \cite{Olver})
 \begin{equation}
 f  (x, z)\sim   {\cal G} (x,z)^{-1/2}  \exp  \Big(-\int_{x_{0}} ^x  {\cal G} (y,z) dy\Big)=: Q  (x, z)
\label{eq:Ans}\end{equation}
as $x\to\infty$.  Here $x_{0}$ is some fixed number and
\[
{\cal G} (x,z)= \sqrt{\frac{b(x)-z}{a(x)} } , \q\Re {\cal G} (x,z)  \geq 0.
\]
 Note that the function $ Q  (x, z)$ (the Ansatz for the  Jost solution $f  (x, z)$) satisfies equation \e{eq:Schr} with a sufficiently good accuracy. Sometimes (if $a(x) \to a_{\infty}$, $b(x)\to 0$ slowly)  it is convenient (see \cite{Y-LR}) to omit the pre-exponential factor $ {\cal G} (x,z)^{-1/2}$ in \e{eq:Ans}.
 
   For $z=\lambda\in {\Bbb R}$, the regular solution $\psi (x,\lambda)$ is a linear combination of the Jost solutions $f(x,\lambda)$ and $\ov{f(x,\lambda)}$ which yields asymptotics of $\psi (x,\lambda)$  as $x\to \infty$.  For example, in the case $a(x)=1$, $b\in L^1 ({\Bbb R}_{+})$ one has
  \begin{equation}
\psi(x,\lambda)\sim   \kappa(\lambda) \sin (\sqrt{\lambda}  x+ \eta (\lambda))
\label{eq:AsReg}\end{equation}
where $\kappa(\lambda)$ and $\eta(\lambda)$  are known as the scattering (or limit) amplitude and phase, respectively.
  If $\Im z\neq 0$, then one additionally constructs, by an explicit formula,  a solution $g(x,z)$ of \e{eq:Schr} exponentially growing as $x\to\infty$. This yields 	asymptotics of  $\psi(x,z)$ for $\Im z\neq 0$.  This scheme was realized in  \cite{Y-LR}.
  
  Let us compare asymptotic formulas \e{eq:A22P+}  and \e{eq:AsReg}. The crucial difference between them is that the phase $\sqrt{\lambda}  x$ in \e{eq:AsReg} depends on $\lambda$ while $\sqrt{|\tau| n}$ in \e{eq:A22P+}   does not depend on $z$.

    \subsection{Scheme of the approach}
 
 An analogy between the equations \e{eq:Jy} and \e{eq:Schr} is of course very well known. However it seems to be never consistently exploited before. In particular, the papers cited above rely on specific methods of difference  equations. Some of these methods are quite ingenious, but, in the author's opinion, the standard approach of differential equations works perfectly well and allows one to study asymptotic behavior of orthogonal polynomials in a very direct way. This approach was already used in the cases of coefficients satisfying $a_{n} \to a_{\infty}>0$, $b_{n} \to 0$ in \cite{Y/LD, Y-LR} and of coefficients satisfying conditions  \e{eq:Gr} with $|\gamma |\neq 1$ and \e{eq:nc} in \cite{nCarl}.

We are applying the same scheme in the  critical  singular case  when conditions  \e{eq:ASa}  and \e{eq:ASb}  are satisfied with $| \gamma |=1$ and $\sigma>3/2$.   Let us briefly describe the main steps of our approach. 

A. 
First, we forget about the orthogonal polynomials $P_{n} (z)$
and distinguish solutions (the Jost solutions) $f_{n}  (z)$ of the  difference  equation \e{eq:Jy} 
 by their asymptotics as $n\to\infty$. This requires a construction of an Ansatz $Q_{n} $   for   the Jost solutions.

 
  B. 
Under assumption \e{eq:nc} this construction (see Sect.~4) is very explicit and, in particular, does not depend on $z\in{\Bbb C}$. In  the case  $\tau<0$, we set
  \begin{equation}
Q_{n}= 
 \nu^n n^{-\sigma /2  + 1/4}     e^{-2 i \sqrt{| \tau | n}  }.
\label{eq:Ans1}\end{equation}
 In  the case  $ \tau >0$, the Ansatz equals 
 \begin{equation}
Q_{n}= \nu^n n^{-\sigma /2  + 1/4}     e^{-2  \sqrt{ \tau  n}  }.
\label{eq:Ans2}\end{equation}
  In both cases
  the relative remainder
\begin{equation}
 r_{n} (z)  : =( \sqrt{ a_{n-1} a_{n}} Q_{n}  )^{-1} \big(a_{n-1} Q_{n-1}   + (b_{n}-z)Q_{n}   + a_{n} Q_{n+1}  \big), \q n\in{\Bbb Z}_{+},
\label{eq:Grr}\end{equation}
belongs to $ \ell^1 ({\Bbb Z}_{+})$.   
  At an intuitive level, the fact that the Ans\"atzen \e{eq:Ans1} and \e{eq:Ans2} do not depend on $z\in{\Bbb C}$ can be explained by the fast growth of the coefficients $a_{n}$ which makes the spectral parameter $z$  negligible in \e{eq:Grr}.
  
  Actually,  the Ans\"atzen we use  are only distantly similar to the Liouville-Green Ansatz    \e{eq:Ans}.  
  On the other hand, for integer $\sigma$, \e{eq:Ans1}  and \e{eq:Ans2} are close to formulas of the Birkhoff-Adams method significantly polished in \cite{W-L} (see also Theorem~8.36 in the book \cite{El}).

C.
Then we make   a multiplicative change of variables
 \begin{equation}
   f_{n} (z)= Q_{n}  u_{n} (z ) 
      \label{eq:Jost}\end{equation} 
which permits us to reduce   the Jacobi equation \e{eq:Jy} for  $  f_{n} (z)$  to a Volterra ``integral" equation for the sequence $u_{n} (z)$.   
This equation  depends of course on the parameters $a_{n}$, $b_{n}$. In particular, it is     different in the cases $\tau < 0$ and  $ \tau > 0$. However in both cases the Volterra equation for $u_{n} (z)$ is standardly solved    by iterations in Sect.~3 which allows us to prove  that it has a solution such that
$ u_{n} (z)\to 1$ as $n\to\infty$. Then the Jost solutions $ f_{n} (z)$  are defined by formula \e{eq:Jost}.


   
   D. 
   The sequence 
     \begin{equation}
\tilde{f}_{n}( z)=\ov{f_n( \bar{z })}
 \label{eq:AcGG}\end{equation}
also    satisfies  the equation \e{eq:Jy}. In the case $\tau <0$, the solutions  $ f_{n} (z )$ and $   \tilde{f}_{n} (z )$ are linearly independent. Therefore it follows from \e{eq:Ans1}   that all solutions  of the Jacobi equation \e{eq:Jy} have  asymptotic behavior \e{eq:A22P+} with some constants $\kappa_{\pm}$.
   
   In the case $ \tau >0 $,   a solution $ g_{n} (z ) $ of \e{eq:Jy} linearly independent with $ f_{n} (z ) $ can be constructed by an explicit formula
    \begin{equation} 
g_{n} (z ) = f_{n} (z )\sum_{m=n_{0}}^n (a_{m-1} f_{m-1}(z) f_{m}(z))^{-1},\q n\geq n_{0} ,
\label{eq:GEg}\end{equation}
where $n_{0}=n_{0}(z)$ is a sufficiently large number. This solution grows faster than any power of $n$ as $n\to\infty$,
 \[
   g_{n} (  z )=  \frac{ \nu^{n}} {2 \sqrt{\tau} } n ^{-\sigma/2 +1/4}   e^{ 2 \sqrt{\tau  n}} (1+ o(1)) .
\]
Since $g_{n} (  z )$ is linearly independent with $f_{n} (  z )$, the polynomials $P_{n} (z) $  are linear combinations of $ f_{n} (  z )$ and $g_{n} (  z )$ which leads to the formula \e{eq:A2P1}.  

Our plan is the following. A Volterra integral equation for $u_{n} (z)$ is introduced and investigated in Sect.~3; the Jost solutions $f_{n}(z)$ are defined Sect.~4;
asymptotics of the orthogonal  polynomials $P_{n} (z) $  is studied in Sect.~5. The   doubly critical case  is discussed in   Sect.~6.
 






 \section{Difference and Volterra equations}
 
 Here we reduce a construction of the Jost solutions $f_{n}  (z)$ of the Jacobi equation \e{eq:Jy} to a Volterra ``integral" equation which is then solved by iterations. This construction works under very general assumptions with respect to the coefficients $a_{n} $ and $b_{n}$. In particular, conditions \e{eq:ASa} and \e{eq:ASb}
 are not required here.

\subsection{Preliminaries}

Let us consider equation \e{eq:Jy}. Note that the values of $f_{m-1}$ and $f_{m }$ for some $m\in{\Bbb Z}_{+}$ determine the whole sequence $f_{n}$ satisfying the difference equation \e{eq:Jy}.

 Let $f=\{ f_{n} \}_{n=-1}^\infty$ and $g=\{g_{n} \}_{n=-1}^\infty$ be two solutions of equation \e{eq:Jy}. A direct calculation shows that their Wronskian
  \[
W [ f,g ]: = a_{n}  (f_{n}  g_{n+1}-f_{n+1}  g_{n})
\]
does not depend on $n=-1,0, 1,\ldots$. In particular, for $n=-1$ and $n=0$, we have
 \begin{equation}
W [ f,g ]= 2^{-1} (f_{-1}  g_{0}-f_{0}  g_{-1}) \q {\rm and} \q  W [ f,g ] = a_{0}  (f_{0}  g_{ 1}-f_{ 1}  g_{0})
\label{eq:Wr1}\end{equation}
(for definiteness, we put $a_{-1}= 1/2$).
Clearly, the Wronskian $W [ f,g ] =0$ if and only if the solutions $f$ and $g$ are proportional.


It is convenient to introduce a notation
 \[
x_{n}'= x_{n+1}  -x_{n}
\]
for the ``derivative" of a sequence $x_{n}$. Then we have
 \[
(x_{n}^{-1})'= - x_{n}^{-1}x_{n+1}^{-1} x_{n}'
\]
and
 \begin{equation}
(e^{x_{n}})'=  (e^{ x_{n}'}-1) e^{x_{n}}.
\label{eq:dife}\end{equation}
 Note also the Abel summation formula (``integration by parts"):
 \begin{equation}
\sum_{n=N }^ { M} x_{n}  y_{n}' = x_{M}  y_{M+1} - x_{N -1}  y_{N}  -\sum_{n=N } ^{M} x_{n-1}'  y_{n};
\label{eq:Abel}\end{equation}
here $M\geq N\geq 0$ are arbitrary, but we have to set $x_{-1}=0$ so that $x_{-1}'=x_{0}$.


To emphasize the analogy between differential and difference equations, we often  use ``continuous" terminology (Volterra  integral equations, integration by parts, etc.) for sequences labelled by the discrete variable $n$. Below $C$, sometimes with indices,  and $c$ are different positive constants whose precise values are of no importance.

In  constructions below, it suffices to consider the Jacobi  equation \e{eq:Jy} for large $n$ only.

\subsection{Multiplicative change of variables}
   
For construction of $f_{n}  (z)$, we will reformulate the problem    introducing a  sequence
\begin{equation}
 u_{n} (z)=  Q_{n} ^{-1}  f_{n} (z), \q n\in {\Bbb Z}_{+}.
\label{eq:Gs4}\end{equation}
  In our construction, the Ansatz $Q_{n}$ does not depend on $z$. In this section,  we do not make any specific assumptions about  the recurrence coefficients $a_{n}$, $b_{n}$ and  the Ansatz $Q_{n}$ except of course that $Q_{n}  \neq 0$; for definiteness, we set $Q_{-1}=1$.  In proofs, we usually   omit the dependence on $z$ in notation; for example, we write $f_{n}$, $u_{n}$, $r_{n}$.


First, we  derive a difference equation for $ u_{n} (z)$.

\begin{lemma}\label{Gs}
Let the remainder $r_{n} (z)$ be defined by formula \e{eq:Grr}.  Set
\begin{equation}
\Lambda_{n}=\frac{a_{n}}{a_{n-1} }  \frac{Q_{n+1}}{Q_{n-1} } 
\label{eq:GL}\end{equation}
and
\begin{equation}
R_{n} (z) = - \sqrt{\frac{a_{n}}{a_{n-1} } }  \frac{ Q_{n} }{ Q_{n-1}}  r_{n} (z).
\label{eq:GL1}\end{equation}
 Then
 equation  \e{eq:Jy} for a sequence $ f_{n} (z)$ is equivalent to the equation
\begin{equation}
\Lambda_{n}( u_{n+1} (z)- u_{n} (z)) -      ( u_{n} (z)- u_{n-1} (z))=   R_{n} (z) u_{n} (z), \q n\in {\Bbb Z}_{+},
\label{eq:DE}\end{equation}
for  sequence  \e{eq:Gs4}. 
 \end{lemma}

\begin{pf}
Substituting  expression $f_{n}  = Q_{n}  u_{n} $ into  \e{eq:Jy} and using the definition \e{eq:Grr},
we see that
\begin{align*}
( \sqrt{a_{n-1}a_{n}}Q_{n} )^{-1}&\Big(  a_{n-1} f_{n-1} + ( b_{n} -z)f_{n} + a_{n} f_{n+1}\Big)
\\=&
\sqrt{\frac{a_{n-1} }{a_{n}}}\frac{Q_{n-1} }{Q_{n}}  u_{n-1} +\frac{b_{n}-z }{ \sqrt{a_{n-1}a_{n}}}  u_{n} +
\sqrt{\frac{a_{n} }{a_{n-1}}}  \frac{Q_{n+1} }{Q_{n}} u_{n+1} 
\\=&
\sqrt{\frac{a_{n-1} }{a_{n}}}\frac{Q_{n-1} }{Q_{n}} ( u_{n-1}-u_{n})+
\sqrt{\frac{a_{n} }{a_{n-1}}}  \frac{Q_{n+1} }{Q_{n}} (u_{n+1}-u_{n}) +r_{n}  u_{n}
\\=&
\sqrt{\frac{a_{n-1} }{a_{n}}}\frac{Q_{n-1} }{Q_{n}}\Big(( u_{n-1}-u_{n})+ \Lambda_{n}( u_{n+1}  - u_{n} )  - R_{n}u_{n}\Big)
  \end{align*}
  where the coefficients $\Lambda_{n}$ and $R_{n}$ are defined by equalities \e{eq:GL} and \e{eq:GL1}, respectively.
 Therefore    the equations \e{eq:Jy}  and \e{eq:DE} are equivalent.
 \end{pf}

   Our goal here to construct solutions of the equation
\e{eq:DE} such that
  \begin{equation}
\lim_{n\to\infty} u_{n} (z)=   1.
\label{eq:DE1}\end{equation}
 Let us set
 \begin{equation}
X_{n}=       \Lambda_{1}  \cdots \Lambda_{n}
\label{eq:DE2}\end{equation}
and 
\begin{equation}
G_{n,m}=  X_{m-1}  \sum_{p=n }^{m-1}   X_{p}^{-1} , \q   m \geq n+1.
\label{eq:DE3}\end{equation}

The following result  will be proven in the next subsection. Note that we    make assumptions only on products $G_{n,m} R_{m}$ but not on factors $G_{n,m} $ and $  R_{m}$ separately.

\begin{theorem}\label{DE}
Set
\begin{equation}
 h_{m} (z)=\sup _{n\leq m-1} |G_{n,m} R_{m} (z) |
\label{eq:DE4}\end{equation}
and suppose that
\begin{equation}
 \{h_{m}(z)\}\in \ell^1 ({\Bbb Z}_{+}).
\label{eq:DE5}\end{equation}
Then equation \e{eq:DE} has a solution $u_{n} (z)$ satisfying an estimate 
\begin{equation}
 |u _{n}(z)  -1 |\leq e^{H_{n}(z)}-1,\q n\geq 0,
\label{eq:Gp9}\end{equation}
where
 \[
 H_n (z)= \sum_{p=n +1}^\infty  h_p (z) .
\]
In particular, condition \e{eq:DE1} holds. 
 \end{theorem}
 
  \subsection{Volterra equation}
 
 A sequence $u_{n}$ will be constructed as a solution of the Volterra ``integral"   equation
 \begin{equation}
u_{n}  (z)= 1+ \sum_{m=n+1}^\infty G_{n,m}  R_m (z) u_{m} (z).
\label{eq:DE6}\end{equation}
This equation can be solved by successive approximations.

  \begin{lemma}\label{GS3p}
  Let  the assumptions of Theorem~\ref{DE}  be satisfied.    Set $u^{(0)}_n =1$   and 
  \begin{equation}
 u^{(k+1)}_{n} (z)= \sum_{m=n +1}^\infty G_{n,m} R_{m} (z)  u^{(k )}_m (z),\q k\geq 0,
\label{eq:W5}\end{equation}
for all $n\in {\Bbb Z}_{+}$. Then  estimates
 \begin{equation}
| u^{(k )}_{n} (z)  |\leq \frac{H_{n}(z)^k}{k!} ,\q \forall k\in{\Bbb Z}_{+},
\label{eq:W6s}\end{equation}
are true. 
\end{lemma}

 \begin{pf}
  Suppose that \e{eq:W6s} is satisfied for some $k\in{\Bbb Z}_{+}$. We have to check 
 the same estimate (with $k$ replaced by $k+1$ in the right-hand side)  for $ u^{(k+1)}_{n}$.  
  According to definitions \e{eq:DE4} and \e{eq:W5}, it  follows from estimate \e{eq:W6s} that
   \begin{equation}
| u^{(k +1)}_{n} |\leq \frac{1}{k!}  \sum_{m=n +1}^\infty  h_m H_{m}^k.
\label{eq:V7}\end{equation}
Observe that
 \[
H_{m}^{k+1}+ (k+1)   h_{m} H_{m}^k  
\leq
H_{m-1}^{k+1},
\]
and hence, for all $N\in{\Bbb Z}_{+}$,
   \[
 (k+1)  \sum_{m=n +1}^N  h_m   H_{m}^k 
 \leq 
 \sum_{m=n +1}^N  ( H_{m-1}^{k+1}- H_{m}^{k+1})
= H_{n}^{k+1}- H_{N}^{k+1}\leq H_{n}^{k+1}.
 \]
Substituting this bound into   \e{eq:V7}, we obtain estimate \e{eq:W6s} for $u^{(k +1)}_{n}$.
    \end{pf}

Now we are in a position to solve equation \e{eq:DE6} by iterations.

  \begin{theorem}\label{GS3}
 Under the assumptions  of Theorem~\ref{DE} the equation  \e{eq:DE6}
  has a  bounded solution $u_{n} (z)$. This solution satisfies an estimate
\e{eq:Gp9}. 
\end{theorem}

 \begin{pf}  Set
     \begin{equation}
   u_{n} =\sum_{k=0}^\infty u^{(k)}_{n} 
\label{eq:W8}\end{equation}
where $u^{(k)}_{n}$ are defined by recurrence relations \e{eq:W5}.
Estimate \e{eq:W6s} shows that this series is absolutely convergent. Using the Fubini theorem to interchange the order of summations in $m$ and $k$, we see that
\[
   \sum_{m=n+1}^\infty G_{n,m}     R_{m}  u_{m}    =   \sum_{k=0}^\infty\sum_{m=n+1}^\infty G_{n,m}   R_{m}  u_{m}^{(k)} = \sum_{k=0}^\infty  u_{n}^{(k+1)}=-1+\sum_{k=0}^\infty  u_{n}^{(k)}=-1+  u_{n}.
\]
This is equation  \e{eq:DE6} for sequence \e{eq:W8}. Estimate \e{eq:Gp9}  also follows from \e{eq:W6s}, \e{eq:W8}.  
\end{pf} 

 \begin{remark}\label{GS3r}
  A bounded solution $u_{n} (z)$ of \e{eq:DE6} is of course unique. Indeed, suppose  that $\{v_{n}\}
  \in \ell^\infty ({\Bbb Z}_{+})$ satisfies the homogeneous equation  \e{eq:DE6}, that is,
   \[
v_{n}   =  \sum_{m=n+1}^\infty G_{n,m}  R_m   v_{m}  
\]
whence
  \[
| v_{n}  |\leq  \sum_{m=n+1}^\infty h_m   |v_{m}  |.
\]
This estimate implies that
\[
| v_{n}  |\leq    \frac{1}{k!}  \big(\sum_{m=n+1}^\infty h_m\big)^k  \max_{n\in {\Bbb Z}_{+}} \{| v_{n}  |\}, \q \forall k \in {\Bbb Z}_{+}.
\]
It follows that $v_{n}=0$. Note however that we do not use the unicity in our construction.
\end{remark}

 
It turns out that the construction above yields a solution of the difference equation  \e{eq:DE}.
 
  \begin{lemma}\label{GS4}
  Let   $G_{n,m}$ be given by formulas   \e{eq:DE2} and \e{eq:DE3}.
Then  a solution $u_{n}  (z)$ of the  integral  equation   \e{eq:DE6}   satisfies also the difference equation  \e{eq:DE}.
 \end{lemma}
 
  \begin{pf}
  It follows from \e{eq:DE6}    that
  \begin{equation}
 u_{n+1} - u_{n}=   \sum_{m=n+2}^\infty (G_{n+1, m}-G_{n, m})R_m u_{m} -
 G_{n, n+1} R_{n+1} u_{n+1}.
  \label{eq:A17Ma}\end{equation}
  Since according to  \e{eq:DE3}  
\[
G_{n+1, m}-G_{n, m}=-  X_{n}^{-1} X_{m-1} \q \mbox{and}\q
G_{n, n+1}=     1 ,
\]
equality \e{eq:A17Ma} can be rewritten as
\begin{equation}
  u_{n+1} - u_{n} =  - X_{n}^{-1} \sum_{m=n+1}^\infty        X_{m-1}  R_{m} u_{m} .
  \label{eq:A17Mb}\end{equation}
Putting together this equality with the same equality for $n+1$ replaced by $n$, we see that
 \[
 \Lambda_{n} ( u_{n+1} - u_{n}) -   ( u_{n} - u_{n-1})=
 - \Lambda_{n} X_{n}^{-1} \sum_{m=n+1}^\infty        X_{m-1}  R_{m} u_{m} +
 X_{n-1}^{-1} \sum_{m=n}^\infty        X_{m-1}  R_{m} u_{m} .
  \]
Since  $X_{n}=  \Lambda_{n} X_{n-1}$ by \e{eq:DE2}, the right-hand side here equals $ R_{n}u_{n}$, and hence the equation obtained coincides with \e{eq:DE}.  
   \end{pf}
   
   Thus putting together Theorem~\ref{GS3} and Lemma~\ref{GS4}, we conclude the proof of Theorem~\ref{DE}.
   
    \begin{remark}\label{GSy}
    It follows from equality \e{eq:A17Mb} that under the assumptions of Theorem~\ref{DE}, we have an estimate
\begin{equation}
| u_{n}'|  \leq \max_{n\in {\Bbb Z}_{+}}\{ |u_{n}| \} \, |X_{n}|^{-1} \sum_{m=n}^\infty  |      X_{m}  R_{m+1} |.
  \label{eq:A1y}\end{equation} 
 \end{remark}

\subsection{Dependence on the spectral parameter}
   
   
   The results above can be supplemented by the following assertion.
   
    \begin{lemma}\label{DE1}
Let, for some open set $\Omega\subset{\Bbb C}$,   the coefficients   $R_{n} (z)$ be analytic functions of $z\in  \Omega$. Suppose that
the assumptions  of Theorem~\ref{DE} are  satisfied uniformly in $z $ on compact subsets of $z\in\Omega$. Then all functions $u_{n}  (z)$ are also analytic in $z\in \Omega$. Moreover, if $R_{n} (z)$ are continuous up to the boundary of $\Omega$  and the assumptions  of Theorem~\ref{DE} are  satisfied uniformly on $\Omega$, then the same is true for the functions $u_{n}  (z)$.
\end{lemma}

  \begin{pf}
 Observe that if the functions $u_{m}^{(k)} (z)$ in \e{eq:W5} depend analytically (continuously)  on $z$, then the function $u_{n}^{(k+1)} (z)$ is also analytic (continuous). Since  the
  series  \e{eq:W8}  converges uniformly, its sum is also an  analytic (continuous)  function.
     \end{pf}
     
     According to \e{eq:Grr}  and \e{eq:GL1} the remainder $R_{n}(z)$ depends linearly on $z$. In this case it is easy to obtain a bound on $u_{n}  (z)$ for large $lz|$. 
     
       \begin{proposition}\label{Gexp}
       Suppose that sequence  \e{eq:DE4}  satisfies a condition
         \begin{equation}
 h _{n}(z)  \leq {\cal H}_{n}(1+|z|) \q \mbox{where}\q \{{\cal H}_{n}\} \in \ell^1 ({\Bbb Z}_{+} ).
\label{eq:ZZ}\end{equation}
 Then,
 for an arbitrary $\varepsilon>0$ and some constants $C_{n}(\varepsilon) $ $($that   do not depend on $z\in{\Bbb C})$, every function $u_{n} (z)$  satisfies an estimate
  \begin{equation}
 |u _{n}(z)  |\leq C_{n}(\varepsilon) e^{\varepsilon |z|}, \q z\in{\Bbb C}.
\label{eq:ss1}\end{equation}
\end{proposition}
     


  \begin{pf}
According to  inequality \e{eq:Gp9} and condition \e{eq:ZZ} we have an estimate
   \begin{equation}
 |u _m(z)  |\leq e^{H_{m}(z)} \leq e^{\varepsilon_m (1+|z|)} \q \mbox{where}\q 
 \varepsilon_m=\sum_{p=m+1}^\infty{\cal H}_p.
\label{eq:ss2}\end{equation}
According to \e{eq:DE5}
  $ \varepsilon_m\to   0$ as $m\to\infty$.
On the other hand, it follows from equation  \e{eq:DE}  that 
\[
|u_{n-1} (z) | \leq (1 + |\Lambda_{n}| + |R_{n} (z)| ) |u_n (z) | + |\Lambda_{n}| |u_{n+1} (z) |.
\]
Iterating this estimate, we find that
   \begin{align}
|u_{n} (z)| &\leq C_{n} (1+ |z|) (|u_{n+1} (z)|+|u_{n+2} (z)|)\leq\cdots
\nonumber \\
&\leq C_{n,k} (1+ |z|)^k (|u_{n+k} (z)|+|u_{n+k+1} (z)|)
 \label{eq:ss2+}\end{align}
for every $k=1,2, \ldots$. For a given $\varepsilon>0$, choose $k$ such that $2\varepsilon_{n+k }\leq \varepsilon$, 
$2\varepsilon_{n+k+1 }\leq \varepsilon$. Then
putting  estimates \e{eq:ss2} and \e{eq:ss2+} together, we see that 
   \[
 |u _n(z)  |\leq 4 C_{n,k} (1+ |z|)^k e^{\varepsilon |z|/2} .
\]
Since $(1+ |z|)^k\leq c_{k} (\varepsilon) e^{\varepsilon |z|/2}$,
this proves \e{eq:ss1}.
    \end{pf}

Functions $u_{n} (z)$ satisfying estimates  \e{eq:ss1} for all $\varepsilon>0$ are known as functions of minimal exponential type.

  \section{Jost solutions } 
  
   In this   section, we first calculate the remainder \e{eq:Grr} for the Ansatz $Q_{n} $ defined by formulas \e{eq:Ans1} or \e{eq:Ans2}. Then
  we make substitution \e{eq:Jost}  and use Theorem~\ref{DE} to construct an appropriate solution of the corresponding  equation \e{eq:DE}.  This leads to Theorem~\ref{GSSx}.

 \subsection{Ansatz}
  
  Let us apply the results of the previous section to recurrence coefficients $a_{n}$, $b_{n}$ satisfying conditions \e{eq:ASa}, \e{eq:ASb} where $| \gamma |=1$.
  First, we exhibit an Ansatz $Q_{n} $ such that the corresponding remainder \e{eq:Grr} satisfies the condition
   \begin{equation}
r_{n} (z) =O (n^{ -\d}), \q n\to\infty,
\label{eq:rem}\end{equation}
for some $\d >3/2$. We emphasize that this estimate with $\d>1$   used  in the non-critical case in \cite{nCarl}  is not sufficient now.   {\it Until Sect.~6, we always suppose that $\tau\neq 0$.   }  We treat the cases $\tau>0$ and $\tau<0$ parallelly  setting  $\sqrt{ \tau } >0$ if $  \tau >0$  and (for definiteness) $\sqrt{ \tau } = i\sqrt{| \tau |}$ if $  \tau <0$. 

Let us seek $Q_{n}$  in the form 
 \begin{equation}
Q_{n}   = \nu^n  n^s  e^{ - \varphi_{n}  } ,\q n\geq 1  ,\q \nu= \sgn\gamma,
\label{eq:Gr1}\end{equation}
$Q_{0} =1$.  We have to calculate the remainder $r_{n} (z)$  and find an exponent $s$ and a sequence $\varphi_{n} $ such that estimate \e{eq:rem} is satisfied with $\d >3/2$.  Set  $\theta_{n}=\varphi_{n+1}-\varphi_{n}$ and  $\varphi_{1}=0$. Then
 \begin{equation}
\varphi_{n} =\sum_{m=1}^{n-1}\theta_{m}  , \q n\geq 2.
\label{eq:GRR}\end{equation}
 The phases $\theta_{n}$ which we choose below (see \e{eq:Grr9}) obey  the conditions $\theta_n = O(n^{-1/2})$ and $\Re \theta_{n}  \geq 0$.

Put
  \[
 \varkappa_{n}=\sqrt{\frac{a_{n+1}}{a_{n }}}.
 \]
It follows from condition   \e{eq:ASa}  that
\begin{equation}
 \varkappa_{n}  =1+ \frac{\sigma/2} {n}+ O (n^{-2}) ,
\label{eq:BX}\end{equation}
and
\[
 ( a_{n} a_{n-1})^{-1/2} =   n^{-\sigma}  \big(1- (\alpha-\sigma/2) n^{-1}+ O (n^{-2}) \big).
\]
Using also \e{eq:ASb}, we see that  sequence  \e{eq:Gr}  satisfies a relation
\begin{equation}
 \gamma_{n}= \nu \big(1+  (\tau/2 ) n^{-1}+ O (n^{-2}) \big)
\label{eq:BX1}\end{equation}
where $\tau$ is defined by equality
\e{eq:BX3}.
 
Let us now calculate the   remainder \e{eq:Grr}.
We write it as
\[
 r_{n} (z) = \sqrt{ \frac {a_{n-1}} {a_{n}}} \frac {Q_{n-1}} {Q_{n}}
 +
  \sqrt{ \frac {a_{n}} {a_{n-1}}} \frac {Q_{n+1}} {Q_{n}}
   +
 \frac {b_{n} -z} {\sqrt{ a_{n-1}a_{n}}}.
\]
 Since  according to \e{eq:Gr1}  and \e{eq:GRR}
\[
\frac{Q_{n+1}}{Q_{n}}= \nu \big( \frac{ n+1}{n}\big)^s e ^{-\theta_{n}},
\]
easy calculations yield the following assertion.

\begin{lemma}\label{RK}
The relative remainder \e{eq:Grr} can be rewritten as
 \begin{equation}
 r_{n} (z)  =    
  - \nu \varkappa_{n-1}^{-1} \big(\frac{n-1}  {n}\big)^s  e^{\theta_{n-1}}    - \nu \varkappa_{n-1}\big(\frac{n+1}  {n}\big)^s e^{-\theta_{n}}    
+ 2\gamma_{n} - z  ( a_{n} a_{n-1})^{-1/2} .
\label{eq:Grv}\end{equation}
 \end{lemma}

    \subsection{Estimate of the remainder}

  Here we estimate expression \e{eq:Grv}. Using relations 
    \[
  \big(\frac{n+1}  {n}\big)^s=1+  \frac{s} {n}+ O (n^{-2}) ,
\]
\e{eq:BX} and setting  $k= s+\sigma/2$,  we see that 
\begin{equation}
 r_{n}  =  -  \nu (1-  k n^{-1} )   e^{\theta_{n-1}}  
  -  \nu (1+  k n^{-1} )   e^{-\theta_{n} } +2 \gamma_{n}
  - z (a_{n} a_{n-1})^{-1/2}  + O (n^{-2}) 
\label{eq:Grr2}\end{equation}
where $\gamma_{n} $ satisfies \e{eq:BX1}.
Since
       \[
       e^{-\theta_{n}} = 1-\theta_{n}+ 2^{-1}\theta_{n}^2- 6^{-1}\theta_{n}^3 + O (n^{-2}),
       \]
expression \e{eq:Grr2} can be written as 
 \begin{equation}
\nu r_{n}  =  r_{n}^{(1)} + r_{n}^{(2)}+ r_{n}^{(3)}    +\tau n^{-1}  - z \nu (a_{n} a_{n-1})^{-1/2}    +O (n^{-2})
\label{eq:rr}\end{equation}
where 
\begin{equation}
r_{n}^{(1)}= 2k n^{-1} \theta_{n} + (1- k n^{-1}) ( \theta_{n}- \theta_{n-1}) ,
\label{eq:rn1}\end{equation}
\begin{equation}
r_{n}^{(2)}= -\theta_{n}^2 + 2^{-1 } (1- k n^{-1})  ( \theta_{n}^2- \theta_{n-1}^2) 
\label{eq:rn2}\end{equation}
and 
\begin{equation}
r_{n}^{(3)}= 6^{-1}  k n^{-1}\theta_{n}^3 +  6^{-1} (1- k n^{-1}) (\theta_{n}^3 -\theta_{n-1}^3 ) .
\label{eq:rn3}\end{equation}
Our goal is to find the numbers $k$ and $\theta_{n}$  such that expression \e{eq:rr}  satisfies estimate \e{eq:rem}.

Let us first consider the quadratic term \e{eq:rn2}. It should cancel with $\tau n^{-1}$, up to a term $O(n^{-2})$. It is convenient to set
 \begin{equation}
 \theta_{n} = 2 \sqrt{  \tau} \big( (n+1)^{1/2} -n^{1/2} \big) .
 \label{eq:Grr9}\end{equation} 
   Then
  \begin{equation}
 \theta_{n} =  \sqrt{  \tau} n^{-1/2}-4^{-1}\sqrt{  \tau} n^{-3/2} + O (n^{-5/2} )
 \label{eq:Grr9x}\end{equation} 
 and
   \[
r_{n}^{(2)}   + \tau n^{-1}     =  O(n^{-2}).
\]

  According to \e{eq:Grr9} for the linear term \e{eq:rn1}, we  have
  \begin{multline*}
r_{n}^{(1)}= 2k  \sqrt{  \tau}  n^{-3/2}
\\
+ (1- k n^{-1}) \sqrt{  \tau} \big( n^{-1/2} -4^{-1}  n^{-3/2} - (n+1)^{-1/2} +4^{-1}  (n+1)^{-3/2}\big)
\\
+ O(n^{-5/2 })= (2k-1/2)  \sqrt{  \tau}  n^{-3/2} + O(n^{-5/2 }) .
\end{multline*}
The coefficient at $n^{-3/2 } $ is zero if $k=1/4$, that is, 
 \begin{equation}
s= -  \sigma /2 +1/4.
 \label{eq:k1}\end{equation}

Finally, for cubic terms \e{eq:rn3}, we  use that
\[
\theta_{n}^3 -\theta_{n-1}^3 =(\theta_{n}  -\theta_{n-1}) (\theta_{n}^2 +\theta_{n}\theta_{n-1}  
+ \theta_{n-1}^2)= O(n^{-3/2}) O(n^{-1})= O(n^{-5/2}),
\]
whence
$r_{n}^{(3)}     = O(n^{-5/2})$.



Let us summarize these calculations and observe that for our choice \e{eq:Grr9} the phase \e{eq:GRR} equals
\begin{equation}
\varphi_{n} = 2 \sqrt{  \tau}  n^{1/2}   .
 \label{eq:k2a}\end{equation} 
 
    \begin{lemma}\label{Gr}
Let conditions  \e{eq:ASa} and \e{eq:ASb}  where $|\gamma|=1$  be satisfied. Define the numbers $\tau$ and $s$ by equalities \e{eq:BX3}  and \e{eq:k1} and suppose that  $\tau\neq 0$. Define the sequence $Q_{n}$ by formula \e{eq:Gr1}   where
$\varphi_{n}     $
is given by \e{eq:k2a}.
Then  the corresponding remainder \e{eq:Grr} satisfies the estimate \e{eq:rem} with $\d = \min\{\sigma , 2\}$.  
\end{lemma}

 According to \e{eq:GL1} we have
\[
R_{n} (z) =-\nu \sqrt{\frac{a_{n}}{a_{n-1} } } \big(\frac{ n }{ n-1} \big)^s e^{ -\theta_{n-1}}r_{n} (z),
\]
whence by  Lemma~\ref{Gr}
   \begin{equation}
R_{n} (z) =O (n^{-\d}) \q \mbox{where}\q \d= \min\{\sigma , 2\}.
\label{eq:Rem}\end{equation}

     \subsection{Estimate of the ``integral" kernel}
     
     Here  we estimate matrix elements $G_{n,m}  $ defined by formulas \e{eq:DE2} and \e{eq:DE3} where $\Lambda_{n}$ is given by \e{eq:GL}.  Now the product \e{eq:DE2}  equals 
\[
X_{n}=      \frac{a_{1}a_{2} \cdots a_{n} } {a_{0}a_{1} \cdots a_{n-1} } 
\;
 \frac{Q_2 Q_3 \cdots Q_{n} Q_{n+1}  } {Q_0 Q_1 \cdots Q_{n-2} Q_{n-1}  }= \frac{a_{n}Q_n  Q_{n+1}  } {a_{0}Q_0 Q_1    } .
\]
It  follows from definition \e{eq:Gr1} that
\begin{align}
X_{n}=        a_{0}^{-1} x_{n}e^{-\varphi_{n}  - \varphi_{n-1}}
\label{eq:Gq1}\end{align}
 where
\begin{equation} 
x_{n}  =  (n+1)^s n^s  a_{n} =  n^{1/2} (1+ c n^{-1}+ O(n^{-2}))  
\label{eq:Gq3}\end{equation}
according to \e{eq:ASa}  and \e{eq:k1};
  the precise value of the constant $ c $ here is inessential.   

Let us state a necessary estimate on $G_{n,m}  $.

\begin{lemma}\label{GS2+}
Let $X_{n}$ be given by formulas \e{eq:Gq1}, \e{eq:Gq3}.  Then
 for all $m>n\geq 0$, matrix elements \e{eq:DE3}      satisfy an estimate
\begin{equation}
|G_{n,m}|\leq    Cm^{1/2}
\label{eq:AbG}\end{equation}
with a constant  $C$ that does not depend on $n$ and $m$.
 \end{lemma}

 \begin{pf}
 By definition  \e{eq:Gq1}, we have 
 \begin{equation}
 a_{0}^{-1}\sum_{p=n}^{m-1}   X_{p}^{-1} =    \sum_{p=n}^{m-1}  x_{p}^{-1}e^{\varphi_{p}  + \varphi_{p-1}}  
\label{eq:Gq4}\end{equation}
where according to  \e{eq:dife}
 \begin{equation}
e^{\varphi_{p}  + \varphi_{p-1}}  = (e^{ \theta_{p}  + \theta_{p-1}} -1)^{-1} (e^{ \varphi_{p}  + \varphi_{p-1}}  )' .
\label{eq:Gn}\end{equation}
Set
  \[
y_{n} = x_n^{-1}  (e^{ \theta_{n}  + \theta_{n-1}} -1)^{-1}.
\]
Using formula \e{eq:Abel}, we can integrate by parts in \e{eq:Gq4} which yields
  \begin{equation}
 a_{0}^{-1}  \sum_{p=n}^{m-1}   X_{p}^{-1} =  y_{m-1} e^{\varphi_{m} +  \varphi_{m-1}}-  y_{n-1} e^{\varphi_{n} + \varphi_{n-1}}
 -   \sum_{p=n}^{m-1}  y_{p-1}' e^{\varphi_{p}  + \varphi_{p-1}}  
\label{eq:Gq6}\end{equation}
  
  Let us estimate the right-hand side of \e{eq:Gq6}.
It follows from formula \e{eq:Grr9x}
 that 
 \begin{equation}
 e^{\theta_{n}  + \theta_{n-1}} -1= 2 \sqrt{\tau} n^{-1/2}  (1+ c_{1} n^{-1/2}+ c_{2} n^{-1}+ O(n^{-3/2}))
\label{eq:Gq8}\end{equation}
where  precise values of the constants $ c_{1}  , c_{2} $ are inessential. In particular, \e{eq:Gq8} implies that 
  \[
  |e^{\theta_{n}  + \theta_{n-1}} -1|  \geq c  n^{-1/2},  \q c>0.
 \]
 Putting together relations \e{eq:Gq3} and \e{eq:Gq8}, we find that 
 \[
y_{n}=(2\sqrt{\tau} )^{-1} \big(1+ d_{1} n^{-1/2}+ d_{2} n^{-1}+ O(n^{-3/2})\big)
\]
for some constants $ d_{1}    $ and $d_{2}$,
 whence
\begin{equation}
| y_{n} |\leq C<\infty \q \mbox{and}  \q y_{n}'=  O(n^{-3/2}).
\label{eq:Gq6a}\end{equation}

Let us multiply equality  \e{eq:Gq6}  by $S_{m-1}$ and take into account that  
  \begin{equation}
| e^{\varphi_{p}  }| \leq |e^{  \varphi_{m}}| , \q p\leq m,
\label{eq:Gx}\end{equation}
because $\Re \theta_{n}  \geq 0$. By definition \e{eq:DE3}, we now have
   \[
|G_{n,m}| \leq |x_{m-1}| \Big( |y_{m-1} | + | y_{n-1} | 
 +   \sum_{p=n}^{m-1}  |y_{p-1}' |  \Big) .
\]
According to   \e{eq:Gq3} the first factor here
 is estimated by $C m^{1/2}$, and according to   \e{eq:Gq6a} the second factor  
 is uniformly bounded.  This yields \e{eq:AbG}.
  \end{pf}

   \subsection{Jost solutions}
   
   Let us come back to the equation \e{eq:DE} with $\Lambda_{n}$ and $R_{n} (z) $  given by \e{eq:GL} and \e{eq:GL1}.


\begin{theorem}\label{GSSy}
       Let the assumptions  \e{eq:ASa},  \e{eq:ASb} with $\sigma>3/2$ and $|\gamma |=1$  be satisfied.  Set
       $\varrho=\min\{\sigma-3/2, 1/2\}$.  Suppose that $\tau\neq 0$.
     For all   $z\in  \Bbb C$,    the equation \e{eq:DE}   has a solution $u_{n}( z )$ with asymptotics
    \begin{equation}
u_{n}(  z )  =   1 + O(n^{-\varrho})   , \q n\to \infty.
\label{eq:Auu}\end{equation} 
Moreover, 
\begin{equation}
u_{n}'(  z )  =     O(n^{-1/2 -\varrho})   , \q n\to \infty.
\label{eq:Avv}\end{equation} 
For all $n\in {\Bbb Z}_{+}$, the functions $u_{n}( z )$ are entire functions of $z\in  \Bbb C$ of minimal exponential type.
 \end{theorem}
 
 \begin{pf}
 Let us proceed from Theorem~\ref{DE}.
   Estimates \e{eq:Rem} and \e{eq:AbG}   
   show the sequence \e{eq:DE4} satisfies a  bound $h_{m} =O (n^{-\sigma+1/2})$. It  is in $\ell^1 ({\Bbb Z}_{+})$ if $\sigma> 3/2$.    Therefore Theorem~\ref{DE} yields a solution $u_{n}$ satisfying estimate \e{eq:Gp9}
  which implies  \e{eq:Auu}.
   
   For estimates of derivatives $u_{n}' $, we use Remark~\ref{GSy}. By virtue of \e{eq:Gq1}, \e{eq:Gq3} and \e{eq:Gx} inequality \e{eq:A1y} implies that
      \[
| u_{n}'|  \leq  C  n^{-1/2} \sum_{m=n}^\infty  m^{1/2}  | R_{m+1} |.
  \]
  Using also  \e{eq:Rem}, we obtain estimate  \e{eq:Avv}.
 \end{pf}
 
 Now we set
   \[
   f_{n}(z)= Q_{n}   u_{n}  (z).
  \]
    According to Lemma~\ref{Gs} this sequence satisfies the Jacobi equation \e{eq:Jy}  and according to definition \e{eq:Gr1} of $Q_{n}$ it has asymptotics \e{eq:A22G+} or \e{eq:A22G-}.  Therefore Theorem~\ref{GSSx}  is a direct  consequence of Theorem~\ref{GSSy}.  The asymptotic relations \e{eq:A22G+} and \e{eq:A22G-} mean that
 the solution $f_{n}  (z)$ is oscillating for $\tau<0$  and tends to zero faster than any power of $n^{-1}$ for $\tau>0$. 




  In the case $\tau>0$, the condition $f_{n}= Q_{n} (1+ o(1))$ as $n\to\infty$ determines a solution of the Jacobi equation \e{eq:Jy} uniquely. Indeed, suppose that  
    \begin{equation}
    f_{n}= Q_{n} u_{n} \q\mbox{and}  \q      \ti{f}_{n}= Q_{n} \ti{u}_{n} 
     \label{eq:WW}\end{equation} 
      where $u_{n}=  1+ o(1)$ and $\ti{u}_{n}=  1+ o(1)$.
   It  follows that
  \begin{equation}
  W[f, \ti{f}]= \lim_{n\to\infty}\big( a_{n}  Q_{n} Q_{n+1} (u_{n}\ti{u}_{n+1}-u_{n+1}\ti{u}_{n})\big) 
  \label{eq:WW1}\end{equation} 
  equals zero
  because according to  \e{eq:Gr1} and \e{eq:k2a} $Q_{n}\to 0$ faster than any power of $n^{-1}$.
  Therefore $\ti{f} _{n}=c f_{n}$ where $c=1$ by virtue of \e{eq:WW}.
  
  Essentially similar result is true  in the case $\tau<0$. Now we have to suppose that $u_{n}$ and $\ti{u}_{n}$ in \e{eq:WW} satisfy conditions \e{eq:Auu} and \e{eq:Avv}. Then
     \begin{equation}
 u_{n}\ti{u}_{n+1}-u_{n+1}\ti{u}_{n} =   u_{n}\ti{u}_{n}' -u_{n}'\ti{u}_{n}=O(n^{-1/2 -\varrho}).
  \label{eq:WW2}\end{equation} 
  According to \e{eq:Gr1}  and \e{eq:k1} we have
    \begin{equation}
   a_{n}  Q_{n} Q_{n+1} = O (a_{n} n^{2s})=O ( n^{1/2}). 
  \label{eq:WW3}\end{equation} 
  Combining the last two relations, we see that the Wronskian \e{eq:WW1} equals zero whence  again $\ti{f} _{n}= f_{n}$.
 

   \section{Orthogonal polynomials} 
        
  As usual, we suppose that conditions  \e{eq:ASa}  and  \e{eq:ASb} are satisfied with $\sigma>3/2$ and $|\gamma|=1$.
 
      \subsection{ Subcritical case  }

We first consider the case $\tau< 0$.  In addition to the Jost solution $f(z)=\{f_n (z)\}$  constructed in Theorem~\ref{GSSx}, 
we can   define the conjugate Jost solution $\tilde{f}_{n}( z)$ by   formula \e{eq:AcGG}.
 It also satisfies equation \e{eq:Jy}  since the coefficients $a_{n}$ and $b_{n}$ are real, and it has the asymptotics
   \[
\tilde{f}_{n}(  z )= \nu^n n^{-\sigma /2  + 1/4} e^{  2 i \sqrt{| \tau | n}} \big(1 + O( n^{-\varrho})\big), \q n\to\infty,
\]

 \begin{lemma}\label{Wro}
 The Wronskian  of the solutions $f(z)$ and $\tilde{f}_{n}( z)=\ov{f_n( \bar{z })}$ of the Jacobi equation  \e{eq:Jy} equals
  \begin{equation}
W[ f(z), \tilde{f} (  z ) ] =    2i  \nu \sqrt{| \tau |}\neq 0
   \label{eq:A2CG1}  \end{equation}
   so that these solutions are linearly independent.  
 \end{lemma}

 \begin{pf}
 Let us use notation  \e{eq:Gr1}. 
 We now have $f_{n}  = Q_{n}  u_{n}$ and
 $\tilde{f}_{n}= \bar{Q}_{n} \tilde{u}_{n} $ where the sequences $u_{n}$ and $\tilde{u}_{n} $ satisfy conditions \e{eq:Auu} and \e{eq:Avv} so   that
     \begin{multline}
  f_{n}\ti{f}_{n+1}-f_{n+1}\ti{f}_{n}= Q_{n}\bar{Q}_{n+1}  u_{n}\ti{u}_{n+1}
  - Q_{n+1}\bar{Q}_{n}  u_{n+1}\ti{u}_{n}
  \\
  = (Q_{n}\bar{Q}_{n+1} -Q_{n+1}\bar{Q}_{n})  u_{n}\ti{u}_{n+1}
  - Q_{n+1}\bar{Q}_{n}  ( u_{n}\ti{u}_{n+1}-u_{n+1}\ti{u}_{n}).
\label{eq:UV}\end{multline} 
Using definitions  \e{eq:Gr1} and \e{eq:k2a}, we see that
\begin{multline}
Q_{n}\bar{Q}_{n+1} -Q_{n+1}\bar{Q}_{n} =-2i n^s  (n+1)^s \sin\big(2\sqrt{| \tau|}(\sqrt{n+1}- \sqrt{n})\big)
\\
= 2i \nu \sqrt{| \tau|} n^{2s -1/2} \big(1+O (n^{-1})\big).
\label{eq:UV1}\end{multline} 
Observe also that
  \begin{equation}
   Q_{n+1}\bar{Q}_{n}  ( u_{n}\ti{u}_{n+1}-u_{n+1}\ti{u}_{n})= O (n^{2s -1/2 -\varrho})
  \label{eq:UV2}  \end{equation}
 according to \e{eq:WW2} and \e{eq:WW3}.  Substituting \e{eq:UV1} and \e{eq:UV2} into \e{eq:UV} and taking into account condition \e{eq:k1}, we arrive at \e{eq:A2CG1}.
   \end{pf}
 
    
   
  Since $\sigma>3/2$, both sequences $f(z)$ and $\tilde{f} (  z )$ belong to $\ell^2 ({\Bbb Z}_{+})$.
  It follows that the minimal Jacobi operator $J_{0}$ has deficiency indices $(1,1)$ whence  all its self-adjoint extensions have discrete spectra in view of general results of  \cite{Nevan}.  This concludes the proof of part~$1^0$ of Theorem~\ref{S-Adj}.

 
 

By virtue of \e{eq:A2CG1}, an arbitrary solution $F(z)=\{F_{n}(z)\} $ of the Jacobi equation  \e{eq:Jy} is a linear combination of the Jost solutions $f(z)$ and $ \tilde{f} (  z )$, that is 
  \begin{equation}
F_{n}(z)= \kappa_{+} (z)f_{n}(z) + \kappa_{-} (z) \tilde{f}_{n}(  z ),
\label{eq:LC}\end{equation}
where the constants can be expressed via the Wronskians:
  \begin{equation}
  \kappa_{+}(z)=   \frac{\{F(z), \tilde{f} (  z )\}} {2i  \nu\sqrt{2|\tau |}}, \q
   \kappa_{-}(z)= -\frac{\{F(z), f(  z )\} }{2i \nu \sqrt{2|\tau |}} .
\label{eq:LC1}\end{equation}
Thus we arrive at the following assertion.

 \begin{theorem}\label{LC}
       Let the assumptions  \e{eq:ASa},  \e{eq:ASb} with $| \gamma | =1$, $\sigma>3/2$   be satisfied,   and let $\tau < 0$. Choose some $z\in {\Bbb C}$ and put $\varrho=\min\{\sigma-3/2, 1/2\}$.
 Then  an arbitrary solution $F_{n}   $ of the Jacobi equation  \e{eq:Jy} has  asymptotics
  \begin{equation}
F_{n} = \nu^n n^{-\sigma /2  + 1/4}   \big(\kappa_{+} e^{ - 2 i \sqrt{| \tau | n}} + \kappa_{-} e^{  2 i \sqrt{|\tau | n}}\big) \big(1 +  O(n^{-\varrho})\big) , \q n\to\infty,
\label{eq:LC2}\end{equation} 
for constants $\kappa_{\pm } $  defined by \e{eq:LC1}. Conversely, for arbitrary $\kappa_{\pm}\in {\Bbb C}$, there exists a solution $F_{n}$  of the  equation  \e{eq:Jy} with asymptotics \e{eq:LC2}. 
 \end{theorem}

 
 
 Recall that the polynomials   $P(z)=\{P_{n}(z)\}$  are the solutions of the Jacobi equation  \e{eq:Jy} satisfying the conditions  $P_{-1}(z)=0$,  $P_{0}(z)=1$.  Therefore formula  \e{eq:A22P+} of Theorem~\ref{GPS} is a particular case of formula  \e{eq:LC2} 
in Theorem~\ref{LC}.   Moreover,  we now have $\ov{\kappa_{-} (z)}=  \kappa_{+}(\bar{z})$ because $\ov{P_{n}(z) }=  P_{n}(\bar{z})$.  In particular, $\ov{\kappa_{-} (\lambda)}=  \kappa_{+}(\lambda)$  if $\lambda\in{\Bbb R}$.

Theorem~\ref{LC} can be supplemented by the following assertion.

 \begin{proposition}\label{LCx}
 Under the assumptions of Theorem~\ref{LC} suppose that 
 a solution $F_{n}   $ of the  equation  \e{eq:Jy} satisfies a bound
   \begin{equation}
F_{n} = o( n^{-\sigma /2  + 1/4})   
\label{eq:LCx}\end{equation} 
as $n\to\infty$.  Then  $F_{n} =0  $  for all $n\in{\Bbb Z}_{+}$.
 \end{proposition}

    \begin{pf}
        Let us proceed from Theorem~\ref{LC}. Comparing relations \e{eq:LC2} and \e{eq:LCx}  we see that
          \begin{equation}
 \kappa_{+} e^{ - 2 i \sqrt{| \tau | n}} + \kappa_{-} e^{  2 i \sqrt{|\tau | n}} = o(1) , \q n\to\infty .
\label{eq:LCy}\end{equation} 
Let us show that this implies equalities $\kappa_{+}=\kappa_{-}=0$. The modulus of the left-hand side of \e{eq:LCy} is minorated by $| |\kappa_{+} |- | \kappa_{-} | |$ whence $ |\kappa_{+} |= | \kappa_{-} | $.  Let $\kappa_{-}= \kappa_{+}e^{i\theta}$ where $\theta\in [0,2\pi)$ so that
   \begin{equation}
| \kappa_{+} e^{ - 2 i \sqrt{| \tau | n}} + \kappa_{-} e^{  2 i \sqrt{|\tau | n}}\, |= | \kappa_{+} | \,
| e^{  4 i \sqrt{|\tau | n}+i\theta}+ 1|.
\label{eq:LCy1}\end{equation} 
Observe that for an arbitrary sequence $\psi_{n}$ such that $\psi_{n}\to\infty$ and $\psi_{n}' \to 0$ as $n\to\infty$, the set $\{e^{i\psi_{n}}\}$ is dense on the unit circle $\Bbb T$. Indeed, choose an  arc $\Delta\subset\Bbb T$. The points $e^{i\psi_{n}}$  rotate around $\Bbb T$ and they cannot jump over $\Delta$ if $|\psi_{n}'| < |\Delta|$. It follows that $e^{i\psi_{n}}\in \Delta$ for some sufficient large $n$ (actually, for an infinite number of $n$). In particular, the sequence  $e^{i\psi_{n}}$  where $\psi_{n} = 4  \sqrt{|\tau | n}+\theta$ cannot converge to $-1$. Now it follows from 
\e{eq:LCy}  and \e{eq:LCy1} that $\kappa_{+}=\kappa_{-}=0$ whence $F_{n}= 0$ according to 
\e{eq:LC}.         
 \end{pf}


     \subsection{Supercritical case}
       
 Here we consider the case      $\tau > 0$. Now the Jost solutions $\{ f_{n}( z )\}$ and $\{ \tilde{f}_{n}( z )\}$ of equation  \e{eq:Jy} have the same asymptotic behavior  \e{eq:A22G-} as $n\to\infty$. Therefore their Wronskian equals zero, and hence they
 coincide. So,  we have to find another  solution $\{ g_{n}( z )\}$ linearly independent with $\{ f_{n}( z )\}$.
  Choose an arbitrary $z\in{\Bbb C}$.
 Asymptotics   \e{eq:A22G-}   implies that $f_{n}(z)\neq 0$ for sufficiently large $n$, say, $n\geq n_{0}= n_{0}  (z)$.  Let us define $\{ g_{n}( z )\}$
  by the formulas 
    \begin{equation} 
g_{n}(z) = f_{n}(z) G_{n}(z)
\label{eq:GH}\end{equation}
and
   \begin{equation} 
 G_{n} (z)=\sum_{m=n_{0}}^n (a_{m-1}f_{m-1}(z) f_{m}(z))^{-1},\q n \geq n_{0} .
\label{eq:GE+}\end{equation}

First, we recall an elementary assertion of a  general nature.

\begin{theorem}[\cite{nCarl}, Theorem~4.8]\label{GE+}
Suppose that a  sequence $ f (z)=\{f_{n}(z)  \}$  satisfies the Jacobi equation \e{eq:Jy}.
Then   the sequence $ g (z)=\{g_{n}(z)  \}$ defined by formulas  \e{eq:GH} and \e{eq:GE+}
 satisfies the same equation  and
   the Wronskian $W[ f(z), g(z) ]=1$. In particular, 
   the solutions $f(z)$ and $g (z)$ are linearly independent.
  \end{theorem}

It remains to find asymptotics   of the sequence $ g_{n} (z) $ as $n\to\infty$. To that end, we will integrate by parts in \e{eq:GE+}. It  follows from relations \e{eq:Jost},  \e{eq:Gr1} and \e{eq:Gn}  that
\[
 (a_{n-1}f_{n-1} f_{n} )^{-1}= - \big( a_{n-1} (n-1)^s n^s  u_{n-1} u_{n}\big)^{-1}  e^{ \varphi_{n-1}+\varphi_{n}}=- t_{n} \big( e^{ \varphi_{n-1}+\varphi_{n}}     \big)'
\]
where
\begin{equation} 
   t_{n} =   \big( a_{n-1} (n-1)^s n^s \big)^{-1}  \big(e^{\theta_{n-1}+ \theta_{n}}-1)^{-1} \big( u_{n-1} u_{n}\big)^{-1} .
\label{eq:GD1}\end{equation}

\begin{lemma}\label{GD}
Sequence \e{eq:GD1} satisfies relations 
 \[
  t_{n} =  \frac{ 1}{2\sqrt{\tau}}  +  O(n^{-\varrho}),\q    t_{n}' =   O(n^{-1/2-\varrho}).
\]
 \end{lemma}
 
 \begin{pf}
 it suffices to use formula \e{eq:Gq3} for the first factor in the right-hand side of \e{eq:GD1}, formula \e{eq:Gq8} --   for the second factor and apply Theorem~\ref{GSSy} to the third factor.
   \end{pf}

Now formula \e{eq:Abel} of integration by parts yields  a representation for the sequence
 \e{eq:GE+}:
   \begin{multline} 
  G_{n}  =-\sum_{m=n_{0}}^n t_{m} \big(e^{ \varphi_{m-1}+\varphi_{m}} \big)'  = -  t_{n} e^{ \varphi_{n}+\varphi_{n+1}}   
+   t_{n_{0}-1} e^{ \varphi_{n_{0}-1}+\varphi_{n_{0}}} 
+  \wt{G}_{n} 
\label{eq:GEf}\end{multline}
where 
\begin{align}
 \wt{G}_{n} = \sum_{m=n_{0}}^n t_{m-1}'e^{ \varphi_{m-1}+\varphi_{m}} .
\label{eq:GEy1}\end{align}

Let us consider the right-hand side of \e{eq:GEf}.
 It follows from from formula \e{eq:k2a} and Lemma~\ref{GD}   that the first  term  has asymptotics  
     \begin{equation}
t_{n} e^{ \varphi_{n}+\varphi_{n+1}}=\frac{ 1}{2\sqrt{\tau}}   e^{4\sqrt{\tau n}}  (1+ O (n^{-\varrho})), \q n\to\infty. 
\label{eq:gar1}\end{equation}
 The second term   does not depend on $n$. 
  Let us show that the remainder $\wt{G}_{n}  $   is also negligible.

      \begin{lemma}\label{GErX}
Let $\wt{G}_{n} (z)$ be given by formula \e{eq:GEy1}  where $t_{n}$ is defined in \e{eq:GD1}.
 Then
     \begin{equation}
| \wt{G}_{n} (z)| \leq C n^{-\varrho} e^{4 \sqrt{\tau n}}  .
\label{eq:gar}\end{equation}
 \end{lemma}
 
  \begin{pf}
  It follows from  Lemma~\ref{GD}  that
     \begin{equation}
| \wt{G_{n}} | \leq  C     \sum_{m=n_{0}}^n  m^{-1/2 - \varrho} e^{4 \sqrt{ \tau m}}=
C     \sum_{m=n_{0}}^n  p_{m} \big(e^{4 \sqrt{\tau m}}\big)'
\label{eq:GEs1}\end{equation}
where 
     \begin{equation}
 p_{m} =m^{-1/2 - \varrho} \big(e^{4 \sqrt{\tau }(\sqrt{ m+1}- \sqrt{  m} )}-1\big)^{-1}.
\label{eq:GEs2}\end{equation}
Integrating by parts in the right-hand side of \e{eq:GEs1}, we obtain, similarly to   \e{eq:GEf}, that
       \begin{equation}
    \sum_{m=n_{0}}^n  m^{-1/2 - \varrho} e^{4 \sqrt{\tau m}}= p_{n}  e^{4 \sqrt{\tau( n+1)}}
    - p_{n_{0}-1}  e^{4 \sqrt{\tau n_{0}}} -  \sum_{m=n_{0}}^n  p_{m-1}'  e^{4 \sqrt{\tau m}} 
 \label{eq:GEs3}\end{equation}
 According to \e{eq:GEs2} we have
 \[
 p_{m}=O(m^{ - \varrho})\q \mbox{and} \q 
 p_{m}'=O(m^{-1 - \varrho}).
 \]
 Therefore the first term in the right-hand side of \e{eq:GEs3} satisfies estimate
 \e{eq:gar} and the sum is bounded by
 \[
 e^{2 \sqrt{2\tau n}}  \sum_{n_{0}\leq m <n/2}  m^{-1 - \varrho} 
 + e^{4 \sqrt{\tau n}}  \sum_{n/2\leq m \leq n}  m^{-1 - \varrho} .
 \]
  The first sum in the right-hand side is bounded and the second one is $O (n^{-\varrho})$.
   \end{pf}
  
      Let us come back to the representations \e{eq:GH} and    \e{eq:GEf}.
       Putting together  relations   \e{eq:gar1} and \e{eq:gar} and using asymptotics  \e{eq:A22G+}  for $f_{n}  (z)$, we obtain the following result.

 \begin{theorem}\label{GEe}
    Let the assumptions of Theorem~\ref{GSSx} be satisfied, and let $\tau  > 0$. Choose some $z\in {\Bbb C}$.  Then the solution \e{eq:GEg} of the Jacobi equation  \e{eq:Jy} satisfies   the asymptotic relation
     \[
  g_{n} (z) = - \frac{ \nu ^n}{2\sqrt{\tau}} n^{-\sigma /2  + 1/4} e^{2\sqrt{  \tau  n}}  (1+ O(n^{-\varrho})) ,\q n\to\infty.
\]
 \end{theorem}
 
 Since $\{ g_{n} (z)\}\not\in \ell^2 ({\Bbb Z}_{+})$, we can state
 
  \begin{corollary}\label{GEc}
If the assumptions of Theorem~\ref{GSSx} are satisfied and   $\tau > 0$, then the minimal Jacobi operator $J_{0}$ is essentially self-adjoint.
 \end{corollary}

    Set $ P(z)=\{ P_{n}(z)\}_{n=-1}^\infty$, $ f(z)=\{ f_{n}(z)\}_{n=-1}^\infty$ and  
        \begin{equation}
\Omega (z):= W [ P(z), f(z) ] = - 2^{-1}f_{-1}(z)
\label{eq:WRG}\end{equation}
where the first formula \e{eq:Wr1}   has been used. The Wronskian $\Omega (z)$ is also known as the Jost function.  Since $f\in \ell ^2 ({\Bbb Z}_{+})$, we see that
 $\Omega(z)=0$ if and only if $z$ is an eigenvalue of the Jacobi operator $J=\clos J_{0}$. Zeros of $\Omega(z)$  are real  because $J$ is self-adjoint.
 By Theorem~\ref{GE+}, the Wronskian $W [ f(z),g(z)]=1$ so that
  \begin{equation}
P_{n} (z)= \omega(z)f_{n} (z)-\Omega (z)g_{n} (z)
\label{eq:WR+}\end{equation}
with $\omega(z)= W[ P(z),g(z) ]$.  Note that 
$\omega (z)\neq 0$ if $\Omega(z)= 0$. Therefore Theorems~\ref{GSSx}   and \ref{GEe} yield formula \e{eq:GEGE+} where $\kappa (z)=-\Omega(z)$. This concludes the proof of Theorem~\ref{GPS}. 
Moreover,  equality \e{eq:WR+} allows us to supplement it by the following result.

\begin{proposition}\label{GE1}
    Let the assumptions of Theorem~\ref{GSSx} be satisfied, and let $\tau > 0$. If $\Omega(z)=0$, then
 \[
 P_{n}(  z )  =   W[ P(z),g(z)] \:\nu^n n^{-\sigma /2  + 1/4} e^{- 2\sqrt{\tau  n}}  \big(1 + O(n^{-\varrho})  \big).
\]
 \end{proposition}
 
 
 The resolvent  of the self-adjoint  operator $J =\clos J_{0}$ can be constructed by  the standard (cf. Lemma~2.6 in \cite{Y/LD}) formulas. Recall that $e_{n}$, $n\in {\Bbb Z}_{+}$, is the canonical basis in the space $\ell^2 ({\Bbb Z}_{+})$.  


   \begin{proposition}\label{resolvent}
   Under the assumptions of Theorem~\ref{GEe}, the resolvent $(J -z )^{-1}$ of the Jacobi operator  $J$  is given by the equalities  
   \begin{equation}
((J -z )^{-1} e_{n}, e_{m})= \Omega(z)^{-1} P_{n} (z) f_{m}(z),\q \Im z \neq 0, 
\label{eq:RRe}\end{equation}
if $n\leq m$ and $((J -z )^{-1} e_{n}, e_{m})=((J -z )^{-1} e_m, e_n)$.  Here $\Omega(z)$ is the Wronskian \e{eq:WRG}.
 \end{proposition}

In view of Theorem~\ref{GSSx}, $f_{n} (z)$ and, in particular,    $\Omega(z)$ are entire functions of $z\in{\Bbb C}$. This allows us to state

\begin{corollary}\label{rdiscr}
The spectrum of the  operator  $J$  is discrete, and its  eigenvalues $\lambda_{1}, \cdots, \lambda_{k}, \ldots$   are given by the equation  $\Omega(\lambda_{k}) =0$.  The resolvent $(J -z )^{-1} $ is an analytic function  of  $z\in{\Bbb C}$ with poles in the points $\lambda_{1}, \cdots, \lambda_{k}, \ldots$.
 \end{corollary}


 In view of formula \e{eq:WRG} and equation \e{eq:Jy} for $n=0$, the equation for eigenvalues of $J$ can be also written  as
  \[
       (b_{0}-\lambda_{k})f_0(\lambda_{k})+ a_{0} f_1(\lambda_{k})=0.
       \]
       It   follows from representation \e{eq:RRe} for $n=m=0$   that the spectral measure of $J$ is given by the standard formula
         \[
\rho(\{\lambda_{k}\})= 2\frac{f_{0}(\lambda_{k})} {\dot{f}_{-1}(\lambda_{k})}  
\]
where $\dot{f}_{-1}(\lambda )$ is the derivative of  $f_{-1}(\lambda )$ in $\lambda$. Alternatively, we have
       \[
\rho(\{\lambda_{k}\})= \big(\sum_{n=0}^\infty P_{n}(\lambda_{k})^2\big)^{-1}. 
\]
A proof of this relation can be found, for example, in \cite{nCarl}, Sect.~4.4.

   \subsection{Operators with discrete spectrum} 
   
   Although the discreteness of the spectrum of the Jacobi operator $J$ was already verified in
        Corollary~\ref{rdiscr}, by variational technique this result can be obtained under fairly more general assumptions.
        The operator $J$ will now be defined via its quadratic form
          \begin{equation}
J [u,u] = (J_{0} u , u)= \sum_{n=0}^\infty b_{n}|u_{n}|^2 +   \sum_{n=0}^\infty( a_{n-1}u_{n-1} +a_{n}u_{n+1} )\bar{u}_{n}
\label{eq:QF}\end{equation}
where $a_{n}$  and $b_{n}$ are     real numbers and $J_{0}$ is given by matrix \e{eq:ZP+}.  The form $J [u,u]$ is defined on the set $\cal D$,
and it is real.  


           \begin{proposition}\label{discr}
                    Suppose that $b_{n}\to \infty$ $($or $b_{n}\to -\infty)$ and 
that
   \begin{equation}
s_{n}:= | b_{n} | -|a_{n-1}|-  |a_{n }|\to\infty
\label{eq:QF1}\end{equation}
as $n\to\infty$.
Then the form \e{eq:QF} is bounded from below $($from above$)$ and closable. The spectrum of the operator $J$ corresponding to this form is discrete.
 \end{proposition}
           
         \begin{pf}   
           By the Schwarz inequality, we have
            \begin{equation}
     2    |  \sum_{n=0}^\infty a_{n}u_{n+1} \bar{u}_{n}| \leq
            |  \sum_{n=0}^\infty (  | a_{n-1}|+| a_{n}| ) |u_{n}|^2  .
            \label{eq:QF3}\end{equation}
            The sum $\sum a_{n-1}u_{n-1}  \bar{u}_{n}
$ satisfies the same estimate.  Suppose, for example,  that $b_{n}\to + \infty$.  Then \e{eq:QF3} yields an estimate
     \begin{equation}
  J [u,u]\geq \sum_{n=0}^\infty s_{n}  |u_{n}|^2 =( S u,u)
\label{eq:QF2}\end{equation}
where $( S f)_{n}= s_{n}  f_{n}$. 
It now follows from condition   \e{eq:QF1} that the  form   $J [u,u]$ is bounded from below. Since the operator $J_{0}$ is symmetric on $\cal D$, the form $J [u,u]$ is closable and hence  it gives rise to a self-adjoint operator $J$.   The inequality \e{eq:QF2} implies that its spectrum is discrete because    the operator $S $   has discrete spectrum.
       \end{pf}

  Proposition~\ref{discr} holds, in particular,  for  Jacobi operators when $a_{n} >0$. It applies directly to the Friedrichs' extension of the operator $J_{0}$, but its conclusion remains true for all extensions $J$ of  $J_{0}$ because the deficiency indices of $J_{0}$ are finite. Note that condition \e{eq:QF1}  does not guarantee that the operator $J_{0}$ is essentially self-adjoint. 
  Indeed, suppose, for simplicity, that $2 b_{n}=\gamma \sqrt{a_{n-1}  a_{n}}$ where  $|\gamma |>1$.  It is shown in \cite{nCarl} (see Corollary~4.22) that the operator $J_{0}$ is essentially self-adjoint if and only if
  \[
  \sum_{n=0}  a_{n}^{-1}  (|\gamma |+\sqrt{ \gamma ^2-1})^{2n}=\infty.
  \]
  However this series converges if $a_{n}\to\infty$ sufficiently rapidly.
  
  If assumptions \e{eq:ASa},  \e{eq:ASb} are satisfied with $\gamma=1$ and $\tau$ is defined by \e{eq:BX3},  then
  \[
  s_{n}  =\tau n^{\sigma-1}  \big(1+ O(n^{-1}) \big) .
  \]
  Therefore condition \e{eq:QF1}  holds true if $\tau>0$ and  $\sigma>1$.
  This concludes the proof of    Theorem~\ref{S-Adj}.  Moreover, we see that the spectra of all self-adjoint extensions $J$ of the minimal operator $J_{0}$ are discrete for all $\sigma>1$.

Note  that  Proposition~\ref{discr} looks similar to Theorem~8 in \cite{H-L}  where, by some reasons,  it was assumed that $a_{n}  n^{-2}\to\infty$ as $n\to\infty$.

 \begin{example}\label{Lag-m}
 Recall that the spectrum of the self-adjoint Jacobi operator $J$ with the coefficients \e{eq:Lag} is absolutely continuous and coincides with $[0,\infty)$. Suppose now that $a_{n}$ are defined by \e{eq:Lag} for some $p > -1$ but
 $
 b_{n}-2n \to\infty
 $
 as $n\to\infty$.  According to Proposition~\ref{discr} the corresponding Jacobi operator $J$ has discrete spectrum.
 \end{example}

 \section{Doubly critical case}
 
 Here we consider the case where the number $\tau$ defined by \e{eq:BX3} is zero.  We do not treat this problem in its full generality. Instead, we exhibit two classes of Jacobi operators satisfying the condition $ \tau =0$ and admitting rather effective spectral analysis.
 Operators from these classes can be reduced to Jacobi operators with zero diagonal elements.
  Probably, the scheme used in the main part of the paper works also in the doubly critical case with the Ansatz defined   by formulas similar to those of Theorem~8.36 (c) of the book \cite{W-L}  where $\sigma$ is integer.

  \subsection{Dediagonalization of Jacobi operators}
  
  Let us proceed from the well known construction (see, e.g.,  \cite{Domb}) which puts into correspondence to an arbitrary self-adjoint Jacobi operator  $\bf J$ with zero diagonal elements ${\bf b}_{n}$   a couple of Jacobi operators $J^{(\pm)}$. It turns out that under fairly general assumptions on off-diagonal elements ${\bf a}_{n}$  of the operator $\bf J$,   the elements $a_{n}^{(\pm)}$ and $b_{n}^{(\pm)}$ of the operators $J^{(\pm)}$ satisfy the   conditions \e{eq:ASa} and \e{eq:ASb} with $\gamma=1$  and $\tau =0$. 
   
  

Let  $e_{n}$, $n=0,1,\ldots$,  be the canonical basis  in the space  $\ell^2 ({\Bbb Z}_{+})$
and let 
      \begin{equation}
  {\bf J} e_{n}= {\bf a}_{n-1}e_{n-1}+  {\bf a}_{n}   e_{n+1}  
  \label{eq:bf0}\end{equation}
  (as usual we put $e_{-1}=0$). Set $U e_{n}= (-1)^n e_{n}$. Then ${\bf J } U= - {\bf J } U$ so that the operators $\bf J$ and $-\bf J$ are unitarily equivalent. In the case  ${\bf b}_{n}=0$, 
  the   orthogonal polynomials ${\bf P}_{n} (z)$  satisfy the identity ${\bf P}_{n} (-z)  = (-1)^n {\bf P}_{n} (z)$.
  
  
 It follows from the definition \e{eq:bf0}  that the operator ${\bf J}^{2}$ acts    by the formula
       \begin{equation}
       {\bf J}^{2} e_{n}={\bf a}_{n-2}{\bf a}_{n-1}e_{n-2}+ ( {\bf a}_{n-1}^{2}+{\bf a}_{n}^{2}) e_{n}+  {\bf a}_{n}{\bf a}_{n+1}e_{n+2}.
\label{eq:red}\end{equation}
Although  ${\bf J}^{2}$ is not a Jacobi operator, it can be reduced to Jacobi operators $J^{(\pm)}$  on the  subspaces ${\cal H}^{(+)}$ and ${\cal H}^{(-)}$ of  $\ell^2 ({\Bbb Z}_{+})$  spanned by the elements $e_{n}$ with even and odd $n$, respectively.  The elements $e_{n}^{(+)}:=e_{2n}$ and $e_{n}^{(-)}:=e_{2n+1}$ where $n=0,1,\ldots$ are the bases in the spaces ${\cal H}^{(+)}$ and ${\cal H}^{(-)}$,  so that each of the subspaces  ${\cal H}^{(\pm)}$ can be identified with the space $\ell^2 ({\Bbb Z}_{+})$.  According to \e{eq:red} ${\cal H}^{(\pm)}$ are  the  invariant subspaces of the operator $ {\bf J}^{2}$ and 
    \begin{equation}
       {\bf J}^{2} e_{n}^{(\pm)}= a_{n-1}^{(\pm)} e_{n-1}^{(\pm)}+ b_{n}^{(\pm)} e_{n}^{(\pm)}+  a_{n}^{(\pm)}e_{n+1}^{(\pm)}
\label{eq:red1}\end{equation}
where
 \begin{equation}
    a_{n}^{(+)}= {\bf a}_{2n} {\bf a}_{2n+1} \q\mbox{and}\q     b_{n}^{(+)}=  {\bf a}_{2n-1}^{2}+{\bf a}_{2n}^{2} 
\label{eq:red2}\end{equation}
and
 \begin{equation}
  a_{n}^{(-)}= {\bf a}_{2n+1} {\bf a}_{2n+2} \q\mbox{and}\q     b_{n}^{(-)}=  {\bf a}_{2n}^{2}+{\bf a}_{2n+1}^{2} .
\label{eq:red2-}\end{equation}
Note that, formally,
\[
    a_{n}^{(-)}=    a_{n+1/2}^{(+)}, \q    b_{n}^{(-)}=    b_{n+1/2}^{(+)}.
    \]
According to \e{eq:red1} the restriction $J^{(\pm)}={\bf J}^{2}\big|_{{\cal H}^{(\pm)}}$  of the operator $ {\bf J}^{2} $ on the subspace ${\cal H}^{(\pm)}$ is a Jacobi operator with matrix elements \e{eq:red2} or \e{eq:red2-}. In particular, we see that the spectral families of the operators ${\bf J}$ and $J^{(\pm)}$  are linked by the formula
  \[
E_{\bf J} (\lambda , \mu)\big|_{{\cal H}^{(\pm)}}=E_{J^{(\pm)}} (\lambda^2, \mu^{2})
\]
where $(\lambda , \mu)\subset {\Bbb R}_{+}$ is an arbitrary interval.

Thus we are led to the following assertion.

\begin{lemma}\label{RED}
 Suppose that a self-adjoint Jacobi operator $\bf J$ has  zero diagonal elements. Let ${\bf a}_{n}$ be its off-diagonal elements, and let the Jacobi operator $ J^{(\pm)}$ have matrix elements \e{eq:red2} or \e{eq:red2-}. Let ${\bf P}_{n} (z)$ and $P_{n}^{(\pm)} (z)$ be the orthogonal polynomials corresponding to the matrix elements $ {\bf a}_{n}$, $ {\bf b}_{n}=0$ and $ a_{n}^{(\pm)}$, $ b_{n}^{(\pm)}$.  Then 
   \begin{equation}
  P_{n}^{(+)} (z) = {\bf P}_{2n}  (\sqrt{z})
\label{eq:red+}\end{equation}
and
  \begin{equation}
  P_{n}^{(-)} (z) = ({\bf a}_{0} \sqrt{z} )^{-1}
{\bf P}_{2n+1}  (\sqrt{z})  .
\label{eq:red-}\end{equation}
  \end{lemma}
  
\begin{pf}
Let us proceed from the Jacobi equation
  \[
{\bf a}_{n-1}  {\bf P}_{n-1}  (z)  + {\bf a}_{n}  {\bf P}_{n+1}  (z)= z {\bf P}_{n}  (z)
\]
for the polynomials ${\bf P}_{n}(z)$. Putting together this equation with the same equations where $n$ is replaced either by $n-1$ or $n+1$, we find that
  \begin{equation}
{\bf a}_{n-1} \big( {\bf a}_{n-2}  {\bf P}_{n-2}  (z)  + {\bf a}_{n-1}  {\bf P}_{n}  (z) \big) 
+ {\bf a}_{n} \big( {\bf a}_{n}  {\bf P}_{n}  (z)  + {\bf a}_{n+1}  {\bf P}_{n+2}  (z) \big) = z^2 {\bf P}_{n}  (z).
\label{eq:bf2}\end{equation}
Setting here $n=2m$ and using notation \e{eq:red2}, we see that
  \[
    a_{m-1}^{(+)}  {\bf P}_{2m-2}  (z) +     b_{m}^{(+)}  {\bf P}_{2m}  (z) 
    +     a_{m}^{(+)}  {\bf P}_{2m+2}  (z)  = z^2 {\bf P}_{2m}  (z).
\]
This is the equation for the orthogonal polynomials $P_{m}^{(+)} (z^2)$ which proves \e{eq:red+}.

Quite similarly, setting   $n=2m+1$ in \e{eq:bf2} and using notation \e{eq:red2-}, we obtain relation \e{eq:red-}.
 \end{pf}
 
 The following example is classical.  Recall that  the Laguerre polynomials $L_{n}^{(p)}(z) $ where the parameter $p>-1$ are defined by the recurrence coefficients \e{eq:Lag}.

 \begin{example}\label{REDe}
 Let ${\bf a}_{n}=\sqrt{(n+1)/2}$, ${\bf b}_{n}=0$. Then $ {\bf P}_{n}  (z)=: H_{n} (z)$ are the Hermite polynomials. 
    The corresponding coefficients \e{eq:red2} and \e{eq:red2-} are given by the  formula \e{eq:Lag} where $p=-1/2$ and $p=1/2$, respectively.
Lemma~\ref{RED}  implies that
\[
H_{2n}  (z)= L_{n}^{(-1/2)} (z^2)\q \mbox{and}\q H_{2n+1}  (z)=\frac{1}{\sqrt{2}} z L_{n}^{(1/2)} (z^2).
\]
These relations are of course very well known (see, e.g., formulas (10.13.2) and (10.13.3) in the book \cite{BE}).
  \end{example}
  
  
The following observation shows that under fairly general assumptions on   ${\bf a}_{n}$,  the asymptotic behavior of the coefficients of the operators $J^{(\pm)}$ is  doubly degenerate.
  
  \begin{lemma}\label{RED1}
  Suppose that
      \begin{equation}
{\bf a}_{n} = \big(n/2\big)^{\sigma/2} \big(  1 +  \boldsymbol{\alpha} n^{-1}+ O (n^{-2})\big)
    \label{eq:bf4}\end{equation}
    for some $\sigma\geq 0$ and $\boldsymbol{\alpha}\in {\Bbb R}$. Define the coefficients $a_{n}^{(\pm)}$  and $b_{n}^{(\pm)}$  by formulas \e{eq:red2} or \e{eq:red2-}. These coefficients satisfy conditions
    \e{eq:ASa} and \e{eq:ASb} with $\gamma^{(\pm)} =1$, $\sigma^{(\pm)}= \sigma$, $\alpha^{(+)}=\boldsymbol{\alpha}+ \sigma/4$,
    $\beta^{(+)}= \boldsymbol{\alpha} - \sigma/4$ and $\alpha^{(-)}= \boldsymbol{\alpha} + \sigma/2$,
    $\beta^{(-)}=\boldsymbol{\alpha}$. For both signs, we have $\tau^{(\pm)}=0$.
  \end{lemma}
  
    \subsection{ Jacobi operators with zero diagonal elements}
    
    To use the results of the previous subsection, we need some information on  orthogonal polynomials satisfying relation \e{eq:RR}  where 
     $b_{n}=0$. In  the Carleman case  \e{eq:Carl}, this class of orthogonal polynomials
  was investigated  in \cite{Jan-Nab, Apt}. It is convenient to state necessary results in the same form as in Sect.~5.2
    of  \cite{nCarl}.


     \begin{theorem}\label{DD0}
     Suppose that  $b_{n}=0$  for all $n\in{\Bbb Z}_{+}$. Let the Carleman condition  \e{eq:Carl} hold and
      \begin{equation}
\frac{a_{n}}{\sqrt{a_{n-1} a_{n+1}}}-1\in \ell^1 ({\Bbb Z}_{+}) .
\label{eq:Gr6b}\end{equation} 
     Set 
     \[
\theta_{n}= 2^{-1} (a_{n}  a_{n-1})^{-1/2}
\] 
and assume that
 \begin{equation}
\theta_{n}' \in \ell^1 ({\Bbb Z}_{+})
\label{eq:DD2}\end{equation} 
 and
\begin{equation}
\theta_{n}^3 \in \ell^1 ({\Bbb Z}_{+}).
\label{eq:DD1}\end{equation} 
 Let
 \begin{equation}
\varphi_{n}= \sum_{m=0}^{n-1} \theta_m  .
\label{eq:psi}\end{equation} 
   Then:
   
   $1^0$   The equation \e{eq:Jy} where $ \pm \Im z \geq 0 $  has a solution $\{f_{n}( z )\}$ with asymptotics
    \begin{equation}
f_{n}(  z )  =  (\mp i)^n a_{n}^{-1 /2   } e^{  \pm  i z \varphi_{n}} \big(1 + o( 1)\big) , \q n\to \infty.
\label{eq:DD}\end{equation} 
In particular, $\{f_{n} (z)\} \in \ell^2 ({\Bbb Z}_{+})$ for $\Im z \neq 0$.
Asymptotics \e{eq:DD} is  uniform in $z$ from compact subsets    of the half-planes $ \pm \Im z \geq 0 $.
For all $n\in {\Bbb Z}_{+}$, the functions $f_{n}( z )$ depend analytically on $z$ for $\pm \Im z >0$ and are continuous up to the real axis $\Im z =0$.

  $2^0$ 
Let the Wronskian  $\Omega (z)$ be defined by formula \e{eq:WRG}. Then  $\Omega (z)\neq 0$ for all $z $  with $\pm \Im z \geq 0$, and the asymptotic behavior as $n\to\infty$ of  the   orthogonal polynomials $P_{n}(z)$  is given by the relations
\begin{equation}
    P_{n}(z)= \Omega(z) \frac{  (\pm i)^{n+1}  }{ 2\sqrt{a_{n}}  }  e^{  \mp iz \varphi_{n}}  (1+o(1)), \q \pm\Im z>0,    
\label{eq:Gas}\end{equation}
and
  \begin{equation}
 P_{n} (\lambda)= - a_{n}^{-1/2} \Big(   | \Omega ( \lambda+i0)|      \sin (\pi n/2  -\lambda\psi_{n} +\arg \Omega(\lambda+i 0) ) + o(1) \Big) , \q \lambda\in {\Bbb R}.
\label{eq:Sz}\end{equation}

    $3^0$   
 The operator $J=\clos J_{0}$ is self-adjoint, and its resolvent $ (J-z)^{-1}$    is determined by the general formula \e{eq:RRe}. 
The spectrum of  the  operator $J$ is absolutely continuous,  coincides with the whole real axis and the spectral measure is given by the expression
  \[
  d\rho(\lambda)=  \pi^{-1}  | \Omega (\lambda+i0) |^{-2} d\lambda.
  \]
 \end{theorem}
 



     \begin{remark}\label{DD1A}
   Theorem~\ref{DD0} remains essentially true for slowly increasing coefficients $a_{n}$ when condition \e{eq:DD1}  is violated. However, formulas \e{eq:DD}, \e{eq:Gas} and \e{eq:Sz} become more complicated  in this case.    
\end{remark}

      \begin{remark}\label{DD1}
      It is easy to see that under condition \e{eq:ASa}  where $\sigma \in (1/3,1]$ all assumptions
  of Theorem~\ref{DD0} on $a_{n}$ are satisfied. In this case, we have
  \[
\varphi_{n}= \frac{n^{1-\sigma}}{2(1-\sigma)}+ C_{\sigma}+ o(1)\; \mbox{if}\; \sigma<1
\q\mbox{and}\q \varphi_{n}=2^{-1}\log n+ C_1+ o(1)\; \mbox{if}\; \sigma=1
\]
for some constants $C_{\sigma}$.
This allows us to simplify formulas  \e{eq:DD}, \e{eq:Gas} and \e{eq:Sz}. For example,  \e{eq:DD} where $\pm\Im z\geq 0$ yields
 \[
f_{n}(  z )  =  (\mp i)^n n^{-\sigma /2   } \exp \big( \pm i z  \frac{n^{1-\sigma}}{2(1-\sigma)}\big) \big(1 + o( 1)\big) \q \mbox{if}\q \sigma<1,
\]
and
 \[
f_{n}(  z )  =  (\mp i)^n n^{( -1 \pm iz)/2  } \big(1 + o( 1)\big) \q \mbox{if}\q \sigma=1.
\]
Similar simplifications can be made in  
  formulas \e{eq:Gas} and \e{eq:Sz}  for the orthogonal  polynomials.
    \end{remark}
    
    \begin{example}[Stieltjes-Carlitz polynomials]\label{Sti}
 One of  the  families  of Stieltjes-Carlitz polynomials is defined  (see formula (9.3) and (9.4) in the book \cite{Chihara}, Chapter~VI, Sect.~9) by the recurrence coefficients $b_{n}=0$ and 
  \begin{equation}
   a_{n}=k(n+1)\; \mbox{for} \; n \, \mbox{even} \q \mbox{and} \q a_{n}=n+1\; \mbox{for}\; n\; \mbox{odd} 
   \label{eq:St-Carl}\end{equation}
 where a parameter $k>0$. 
     Theorem~\ref{DD0} applies for $k=1$; in this case the spectrum of the self-adjoint Jacobi operator $J=\clos J_{0}$  is absolutely continuous and  coincides with the whole real axis. On the contrary, if $k\neq 1$, then (see \cite{Chihara}, Chapter VI, Sect.~9) the spectrum of the   operator $J $ is discrete and consists of eigenvalues $c(k) (j+1/2)$ if $k<1$ and eigenvalues $c(k) j$ if $k>1$; here $j\in {\Bbb Z} $ and $c(k)$ are some constants. Note that for $k\neq 1$ conditions  \e{eq:Gr6b}  and \e{eq:DD2}   are violated. Example \e{eq:St-Carl} shows that these conditions  cannot be omitted in  Theorem~\ref{DD0}. 
     
     This example   exhibits a curious phenomenon: the absolutely continuous spectrum of $J$ for $k=1$ is transformed  into  a   purely discrete  spectrum by an arbitrary small  perturbation of $k$.
\end{example}

  Formulas \e{eq:Gas} and \e{eq:Sz}  are of course consistent with the classical asymptotic expressions for the Hermite polynomials  when $a_{n }=\sqrt{(n+1)/2}$ and $b_{n}= 0$ (see, e.g., Theorems~8.22.6 and 8.22.7 in the G.~Szeg\H{o}'s book \cite{Sz}).

  The Carleman condition \e{eq:Carl}  is not very important here. Under assumption \e{eq:nc}   the sequence \e{eq:psi} has a finite limit as $n\to\infty$, and hence it can be omitted in asymptotic formulas of Theorem~\ref{DD0}.  On the other hand, spectral results are drastically different in these cases. The following assertion is a particular case of Theorems~3.9 and 4.2 in \cite{nCarl}.
  
   \begin{theorem}\label{DD-nC}
     Let  $b_{n}=0$.
     Suppose that conditions \e{eq:nc} and \e{eq:Gr6b}   are satisfied.  
   Then:
   
   $1^0$   For all  $z\in {\Bbb C}$, the equation \e{eq:Jy}   has  a solutions $\{f_{n} ( z )\}$ with asymptotics
     \[
f_{n} (  z )  =  (- i )^n a_{n}^{-1 /2   }  \big(1 + o( 1)\big) , \q n\to \infty.
\]
The sequence $\tilde{f}_{n} (z)=\ov{f_{n} (\bar{z})}$ also satisfies \e{eq:Jy}.
The functions $f_{n} ( z )$ and $\tilde{f}_{n} ( z )$ are analytic in the whole complex plane $\Bbb C$.  Both solutions 
$\{f_{n}  (z)\}, \{\tilde{f}_{n}  (z)\}$ are in $ \ell^2 ({\Bbb Z}_{+})$.  
 
  $2^0$ 
The asymptotic behavior as $n\to\infty$ of  the   orthogonal polynomials $P_{n}(z)$  is given by the relation
\[
    P_{n}(z)=\frac{ 1}{ \sqrt{a_{n}}  } \big( \kappa_{+}(z)(- i)^{n}  + \kappa_{-} (z) i^{n} \big) (1+o(1))
     \]
 for some constants  $\kappa_{+} (z)\in{\Bbb C}$ and $\kappa_{-} (z) = \ov{\kappa_{+}(\bar{z})}$.

    $3^0$   
 The symmetric operator $J_{0}$ has a one parameter family of self-adjoint  extensions $J$. All operators $J$ have discrete spectra.
   \end{theorem}

      \subsection{Regular and singular cases}
      
      Here we combine the results of Lemma~\ref{RED} and Theorems~\ref{DD0} or \ref{DD-nC}. 
      We do not discuss the Jost solutions and state asymptotic results in terms of the orthogonal polynomials only.  Let us distinguish regular and singular cases. The first result, for the regular case,  follows from Theorem~\ref{DD0}.

       \begin{theorem}\label{DD3}
 Let condition \e{eq:bf4}  be satisfied with some
   $\sigma\in (2/3,2]$.
  Define the coefficients $a_{n}^{(\pm)}$ and $b_{n}^{(\pm)}$ by one of the formulas \e{eq:red2} or \e{eq:red2-}. Let $J^{(\pm)}$ be the Jacobi operators with these coefficients, and let $P_{n}^{(\pm)} (z)$ be the corresponding orthogonal polynomials.  Then:

    $1^0$    If $z\in {\Bbb C}\setminus [0,\infty)$ and $\Im \sqrt{z}>0$, then
       \[
    P_{n}^{(\pm)}(z)= \kappa^{(\pm)} (z)   (-1)^{n}   n^{-\sigma/4} 
         \exp \big( - i \sqrt{z}  \frac{n^{1-\sigma/2}}{2-\sigma}\big) \big(1 + o( 1)\big), \q   \sigma<2,    
\]
and
  \begin{equation}
    P_{n}^{(\pm)}(z)= \kappa^{(\pm)} (z)   (-1)^{n}   n^{-1/2-i\sqrt{z}} 
   \big(1 + o( 1)\big), \q   \sigma=2,    
\label{eq:GasCx}\end{equation}
as $n\to\infty$
for some constants $\kappa^{(\pm)} (z)$.

 $2^0$    If $\lambda\geq 0$, then
       \[
    P_{n}^{(\pm)}(\lambda)= \kappa^{(\pm)} (\lambda)   (-1)^{n}   n^{-\sigma/4} \big(
         \sin \big(  \sqrt{\lambda}  \frac{n^{1-\sigma/2}}{2-\sigma}+\eta^{(\pm)}(\lambda)\big)   + o( 1)\big), \q   \sigma<2,    
\]
and
  \begin{equation}
   P_{n}^{(\pm)}(\lambda)= \kappa^{(\pm)}(\lambda)   (-1)^{n}   n^{-\sigma/4} \big(
         \sin \big(  \sqrt{\lambda}  \log n+\eta^{(\pm)} (\lambda)\big)   + o( 1)\big), \q   \sigma=2,    
\label{eq:SzCx}\end{equation}
as $n\to\infty$
for some constants $c^{(\pm)} (\lambda)$ and $\eta^{(\pm)} (\lambda)$.

  $3^0$ The operators $J_{0}^{(\pm)}$ are essentially self-adjoint on the set $\cal D$.
        The spectra of their closures  are absolutely continuous and   coincide with $[0,\infty)$.
       \end{theorem}
         
          \begin{remark}\label{LL3}
          According to Lemma~\ref{RED1} condition \e{eq:bf4} ensures that the coefficients $a_{n}^{(\pm)}$ and $b_{n}^{(\pm)}$  satisfy relations \e{eq:ASa} and \e{eq:ASb} with $\gamma^{(\pm)}  =1$, $\sigma^{(\pm)} =\sigma$ and $\tau^{(\pm)} = 0$. We emphasize that, in the doubly critical case, the role of the Carleman condition is played by the assumption $\sigma\leq 2$. Asymptotic phases in formulas of Theorem~\ref{DD3} depend on the spectral parameter. So, we are in the regular situation here.
 \end{remark}

  \begin{remark}\label{LL}
  We emphasize that Theorem~\ref{DD3} yields asymptotics as $n\to\infty$ of the orthogonal polynomials $P_{n} (z)$ uniformly in a neighborhood of the point $z=0$. This is a very strong but specific result which cannot be generically true. For example, for the Laguerre polynomials $L_{n}^{(p)}(z)$ it is true for $p=-1/2$ and $p=1/2$ only. For other values of $p$,  the asymptotics of $L_{n}^{(p)}(z)$ in a neighborhood of the point $z=0$ is given by    a more complicated Hilbs formula (see, e.g., formula (10.15.2) in the book \cite{BE}).
\end{remark}

        \begin{example}\label{LL1}
         Let 
          \begin{equation}
         {\bf a}_{2n-1}= {\bf a}_{2n}=\sqrt{n+1}
            \label{eq:LagV1}\end{equation}
         for all $n\in{\Bbb Z}_{+}$.  Then
           \[
    a_{n}^{(+)}= \sqrt{(n+1)(n+2)}\: , \q     b_{n}^{(+)} =  2n+ 2
    \]
    and the corresponding orthogonal polynomials are Laguerre polynomials $L_{n}^{(1)}(z)$.
    Therefore formulas \e{eq:DD}, \e{eq:Gas} and \e{eq:Sz}  are not true in this case for $z=0$. This does not contradict Theorem~\ref{DD0}  because
       condition \e{eq:Gr6b} is violated for the coefficients \e{eq:LagV1}.
\end{example}


\begin{example}\label{LL2}
         Let 
          \begin{equation}
    a_{n}=\sqrt{(n+1)  (n+x) (n+y)(n+x+y) }\q\mbox{and}\q b_{n}=2n^2+ (2x+2y-1)n +xy
            \label{eq:LagV2}\end{equation}
            where parameters $x,y\in{\Bbb R}_{+}$.
            The corresponding orthogonal polynomials are known as continuous dual Hahn polynomials.  Now assumptions  \e{eq:ASa} and  \e{eq:ASb} are satisfied with  $\gamma=1$, $\sigma=2$,  $\alpha=x+y+ 1/2$ and $\beta=x+y-1/2$ so that $\tau=0$ and   we are in the  doubly critical   case. Note that relations \e{eq:red2} hold true for the coefficients $ a_{n}^{(+)}= a_{n} $  and $ b_{n}^{(+)}= b_{n} $ defined by \e{eq:LagV2} if 
 \begin{equation}
    {\bf a}_{2n}=\sqrt{ (n+x) (n+y)}\q\mbox{and}\q     {\bf a}_{2n+1}=\sqrt{ (n+1) (n+x+y)}.
            \label{eq:LagV3}\end{equation}
       Since
            \[
            {\bf a}_{n}= 2^{-1} \, n \: \big(1 + (x+y) n^{-1} +O (n^{-2})    \big) 
            \]
            both for even and odd $n\in {\Bbb Z}_{+}$,
            the conditions of  Theorem~\ref{DD0}  are satisfied for the coefficients \e{eq:LagV3}.
            Therefore asymptotics of the continuous dual Hahn polynomials are given by formulas 
            \e{eq:GasCx} and    \e{eq:SzCx} and hence the regular case occurs. This is consistent with  the condition $\sigma=2$. The Jacobi operators $J_{0} $ with matrix elements \e{eq:LagV2} are  essentially self-adjoint on the set $\cal D$, and their  closures have absolutely continuous  spectra   coinciding with $[0,\infty)$.
\end{example}

In the singular doubly critical case, we have the following result.  It is  a consequence of Theorem~\ref{DD-nC}.

 \begin{theorem}\label{DD3-nC}
 Let condition \e{eq:bf4}  be satisfied with some
   $\sigma>2$.
  Define the coefficients $a_{n}^{(\pm)}$ and $b_{n}^{(\pm)}$ by one of the formulas \e{eq:red2} or \e{eq:red2-}. Let $J^{(\pm)}$ be the Jacobi operators with these coefficients, and let $P_{n}^{(\pm)} (z)$ be the corresponding orthogonal polynomials.  Then:
   
               $1^0$ 
The asymptotic behavior as $n\to\infty$ of  the   orthogonal polynomials $P_{n}^{(\pm)}(z)$  is given by the relation
\[
    P_{n}^{(\pm)}(z)=  \kappa^{(\pm)} (z)  (- 1)^{n} n^{-\sigma/4}    \big(1+o(1)\big)
     \]
 for some constants  $\kappa^{(\pm)}  (z)\in{\Bbb C}$.

        $2^0$ The operators $J_{0}^{(\pm)}$ have deficiency indices $(1,1)$. All their self-adjoint extensions  $J^{(\pm)}$ have discrete spectra.
            \end{theorem}
            
               \subsection{Comments}
               
               Let  us come back to Figure~1 and summarize the results obtained. Suppose that 
               assumptions \e{eq:ASa}, \e{eq:ASb} are satisfied.
               
                In the non-critical situation $|\gamma| \neq 1$, the singular case occurs if $\sigma>1$. 
                The corresponding Jacobi operators $J_{0}$ are essentially self-adjoint for $|\gamma| >1$, and they have deficiency indices $(1,1)$ for $|\gamma| <1$. In both cases the spectra of self-adjoint extensions of $J_{0}$ are discrete.
                
                   In the critical situation where $|\gamma| = 1$ but $\tau\neq 0$, the singular case occurs if $\sigma>3/2$.      The corresponding Jacobi operators $J_{0}$ are essentially self-adjoint for $\tau > 0$, and they have deficiency indices $(1,1)$ for $\tau< 0$. In both cases the spectra of self-adjoint extensions of $J_{0}$ are discrete.      Surprisingly,   this is, in general,  no longer true in the intermediary case $\tau =0$.  
                   
                      In the doubly critical situation where $|\gamma| = 1$ and $\tau= 0$, the singular case occurs if $\sigma>2$.     Under the assumptions of Theorem~\ref{DD3-nC},
                      the corresponding  operators $J_{0}$  have deficiency indices $(1,1)$ and the spectra of their self-adjoint extensions  are discrete.        According to Theorem~\ref{DD3}  we are in the regular case if $\sigma\leq 2$. The  corresponding  operators $J_{0}$ are essentially self-adjoint  and have absolutely  continuous spectra coinciding with the half-axis $[0,\infty)$ if $\gamma=1$ (or with  the half-axis $(-\infty,1]$ if $\gamma=-1$).
                    
                      We  emphasize that in the doubly critical case, the condition $\sigma>2$ plays the role of $\sigma>3/2$ in the (simply) critical case and  of $\sigma>1$ in the non-critical case.



    \end{document}